\documentclass[11pt]{article}
\usepackage{latexsym,amssymb,amsmath,amsfonts,amsthm}
\makeatletter \topmargin=0mm\headheight=0mm\headsep=0mm
\newcommand{\md}{\mbox{d}}
\newcommand{\be}{\begin{eqnarray}}
\newcommand{\ee}{\end{eqnarray}}
\newcommand{\mr}{\mathbb{R}}
\newcommand{\mc}{\mathcal}
\newcommand{\mb}{\mathbb}
\newcommand{\mx}{\mbox}

\newcommand{\tr}{\triangle}

\newcommand{\ds}{\displaystyle}
\newcommand{\no}{\nonumber}
\newcommand{\ben}{\begin{eqnarray*}}
\newcommand{\enn}{\end{eqnarray*}}

\newcommand{\al}{\alpha}
\newcommand{\la}{\lambda}

\newtheorem{thm}{\textbf Theorem}[section]
\newtheorem{lem}{\textbf Lemma}[section]
\newtheorem{rem}{\textbf Remark}[section]

\newtheorem{prop}{\textbf Proposition}[section]
\newtheorem{defin}{\textbf Definition}[section]

\begin{document}
\begin{titlepage}
\title{\bf Well-posedness of the Cauchy problem for the fractional power dissipative equations}
\author{Changxing Miao$^1$, Baoquan Yuan$^2$ and Bo Zhang$^2$\\
       \\
       $^1$ Institute of Applied Physics and Computational Mathematics,\\
        P.O. Box 8009, Beijing 100088, P.R. China.\\
        (miao\_{}changxing@iapcm.ac.cn)\\
        \\ $^2$ Institute of Applied Mathematics,\\Academy of Mathematics \& Systems Science,\\
          Chinese Academy of Sciences, Beijing 100080, P.R. China\\
          (bqyuan@hpu.edu.cn  and b.zhang@amt.ac.cn)}

\date{}
\end{titlepage}
\maketitle

\begin{abstract}
This paper studies the Cauchy problem for the nonlinear fractional
power dissipative equation $u_t+(-\triangle)^\alpha u= F(u)$ for
initial data in the Lebesgue space $L^r(\mr^n)$ with $\ds r\ge
r_d\triangleq{nb}/({2\alpha-d})$ or the homogeneous Besov space
$\ds\dot{B}^{-\sigma}_{p,\infty}(\mr^n)$ with
$\ds\sigma=(2\alpha-d)/b-n/p$ and $1\le p\le \infty$, where
$\alpha>0$, $F(u)=f(u)$ or $Q(D)f(u)$ with $Q(D)$ being a
homogeneous pseudo-differential operator of order
$d\in[0,2\alpha)$ and $f(u)$ is a function of $u$ which behaves
like $|u|^bu$ with $b>0$.
\vskip0.1in


\vspace{.1in}

\end{abstract}


\section{Introduction}
\setcounter{equation}{0}

In this paper we study the Cauchy problem for the semi-linear
fractional power dissipative equation
\begin{eqnarray}\label{1.1}
\begin{cases}
u_t+(-\triangle)^\alpha u= F(u),\;\; (t,x)\in\mr^+\times\mr^n,\\
u(0,x)=\varphi(x),\;\; x\in \mr^n
\end{cases}
\end{eqnarray}
with the nonlinear term $F(u)$ is equal either to $f(u)$ or to $Q(D)f(u)$, where $Q(D)$ is
a homogeneous pseudo-differential operator of order $d\in[0,2\alpha)$ with real number $\alpha>0$
and $f(u)$ is a function of $u$ which behaves like $|u|^bu$ or $|u|^{b_1}u+|u|^{b_2}u$
with $b>0$, $b_1>0$ and $b_2>0$.
The evolution equation in (\ref{1.1}) models several  classical equations, for example,

(1) the semi-linear fractional power dissipative equation
\ben
u_t+(-\triangle)^\alpha u=\pm\nu |u|^bu.
\enn

(2) the dissipative quasi-geostrophic (QG) equation
\be\label{1.3}
\begin{cases}
\theta_t+u\cdot\nabla\theta+\kappa(-\triangle)^\alpha
\theta=0,\\
u=(u_1,u_2)=\nabla^\bot\psi,\quad (-\tr)^{1/2}\psi=\theta,
\end{cases}
(t,x)\in \mr^+\times\mr^2,
\ee
where $1/2<\alpha\le 1$.

(3) the generalized Navier-Stokes equation
\ben
u_t+(-\triangle)^\alpha u-(u\cdot\nabla)u+\nabla P=0,\quad \nabla u=0.
\enn

(4) the generalized convection-diffusion equation
\ben
u_t+(-\triangle)^\alpha u=a\cdot\nabla(|u|^bu),\
a\in\mr^n\mx{/}\{0\}.
\enn

(5) the Ginzburg-Landau equation
\ben
u_t+a_1\nabla^4 u=G(u)+a_2\nabla^2u+a\nabla^2 u^3,\;\;(t,x)\in [0,\infty)\times\mr^n,
\enn
where $a_1>0$, $a>0$ and $a_2\not=0$.

The case $\alpha=1$ for the problem (\ref{1.1}) corresponds to the classical semi-linear heat
equation and has been studied extensively (see, e.g.
\cite{Fujita}-\cite{M1},\cite{M-Y,M-Z,Ponce,Ribaud1,Ribaud2,Terraneo,W}).
Concerning the generalized Navier-Stokes equations, please refer to \cite{C-K}.
In the case when $\alpha$ is an integer, \cite{Miao1,M-G} established the space-time estimates and
well-posedness of strong solutions in Lebesgue spaces to the problem (\ref{1.1}).
For general $\alpha$,  \cite{M-Y} studied the global well-posedness of solutions to
(\ref{1.1}) for small initial data in pseudomeasure sapces.
In the case $1/2<\al\le 1$ (i.e. in the case of dissipative quasi-geostrophic equation (\ref{1.3})),
well-poswedness of solutions has been studied in, e.g. Lebesgue spaces \cite{Wu1}, Sobolev spaces
\cite{Wu3}, H\"older spaces \cite{Wu5}, Besov spaces \cite{Wu4,Zhang} and Trieble spaces \cite{Chae1}.

In this paper we shall give a unified method to deal with the well-posednesss of
the Cauchy problem (\ref{1.1}) for initial data in the Lebesgue space $L^r(\mr^n)$
($r\ge r_0\triangleq{nb}/({2\alpha-d})$) or in the Besov space $\dot{B}^{-\sigma}_{p,\infty}(\mr^n)$
($\sigma=({2\alpha-d})/b-n/p$ and $1\le p\le \infty$), employing appropriate space-time spaces such as
$C([0,\infty);L^r(\mr^n))\cap L^q([0,\infty);L^p(\mr^n))$ or
$C(I;L^r(\mr^n))\cap \mc{C}_q(I;L^p(\mr^n))$.
In Section 2, we give a detailed analysis of the kernel function of the fractional
power operator semigroup $S_\al(t)=\mx{e}^{t(-\tr)^\alpha}$. In particular,
we derive the point-wise estimates of the kernel function of the semigroup $S_\al(t)$
by an invariant derivative technique (see Lemma \ref{lem2.1} below) which leads to an
equivalent characterization of the Besov space (see \cite{Zhang} for the special case
$1/2<\alpha\le 1$ and $n=2$).
In Section 3 making use of the point-wise estimates of the kernel function obtained in Section 2
we establish the space-time estimates for the corresponding linear fractional power dissipative equation.
Section 4 is devoted to the well-posedness in Lebesgue spaces of the Cauchy problem (\ref{1.1}),
using the space-time estimates established in Section 3 in conjunction with the Banach
contraction mapping principle. In Section 5, we consider the fractional power dissipative equations
with more general nonlinear terms. In particular, the interaction between two different nonlinear terms
is discussed, and the local and small global well-posedness of solutions are established.
Finally, in Section 6 we establish the well-posedness of solutions to the fractional power dissipative
equation (\ref{1.1}) for initial data in the Besov space $\dot{B}^{-\sigma}_{p,\infty}(\mr^n)$
or in the critical Lebesgue space $\ds L^{\frac{nb}{2\alpha-d}}(\mr^n)$ but with
small norm in the Besov space $\dot{B}^{-\sigma}_{p,\infty}(\mr^n)$). Since the Besov space
$\dot{B}^{-\sigma}_{p,\infty}(\mr^n)$ contains self-similar initial data in the sense that the initial
data $\varphi(x)$ satisfies $\ds\lambda^{\frac{2\al}b}\varphi(\la x)=\varphi(x)$ for any $\lambda>0$,
then our results in Section 6 implies the existence of global self-similar solutions to
(\ref{1.1}). Concerning the systematic scaling analysis of nonlinear parabolic equations
please refer to Karch \cite{Karch1,Karch2}.

%


\section{Analysis of the operator semigroup $S_{\alpha}(t)$}
\setcounter{equation}{0}

In this section we consider the linear semigroup $ S_{\alpha}(t)\triangleq\mx{e}^{-t(-\tr)^\alpha}$
generated by the following linear fractional power dissipative equation (\ref{2.1}).
We show that the kernel function of the operator semigroup $S_{\alpha}(t)$
generates a bounded linear operator on $L^p(\mr^n)$ for $p\in [1,\infty]$.

Consider the Cauchy problem for the homogeneous linear fractional power dissipative equation
\begin{eqnarray}\label{2.1}
\begin{cases}
u_t+(-\tr)^\alpha u=0,\ &(t,x)\in [0,\infty)\times\mr^n;\\
u(0)=\varphi(x),\ & x\in \mr^n.
\end{cases}
\end{eqnarray}
By the Fourier transform the solution of the problem (\ref{2.1}) can be written as
\be\label{2.2}
u(t,x)=\mc{F}^{-1}\bigg(\mx{e}^{-t|\xi|^{2\alpha}}\mc{F}\varphi(\xi)\bigg)
=\mc{F}^{-1}\mx{e}^{-t|\xi|^{2\alpha}}*\varphi(x)
\triangleq K_t(x)*\varphi(x).
\ee
Here $\mc{F}$ and $\mc{F}^{-1}$ denote the Fourier and inverse
Fourier transforms, respectively, defined by
\ben
\mc{F}(f)&=&\hat{f}(\xi)=\frac1{(2\pi)^{n/2}}\int_{\mr^n}\mx{e}^{-ix\cdot\xi}f(x)\md x,\\
\mc{F}^{-1}(g)&=&\check{g}(\xi)=\frac1{(2\pi)^{n/2}}\int_{\mr^n}\mx{e}^{ix\cdot\xi}g(\xi)\md\xi
\enn
for any $f,\,g\in\mc{S}'$, where $\mc{S}'$ denotes the space of tempered distributions.

It is well known that for $\alpha=1$ and $\alpha=\frac12$,
$K_t(x)$ is the Gaussian and Poisson kernel function,
respectively, and their properties have been fully understood. In
what follows we consider the general case $\alpha\in (0,\infty)$.
From (\ref{2.2}) and Young's inequality it is seen that, to
guarantee the $L^p\rightarrow L^p$ boundedness of the linear
operator $S_{\alpha}(t)$ one needs only that the kernel
function $K_t(x)$ is bounded on $L^1$. By scaling we have
\be\no
K_t(x)&=&(2\pi)^{-n/2}\int_{\mr^n}\mx{e}^{ix\cdot\xi}\mx{e}^{-t|\xi|^{2\alpha}}\md\xi\\ \no
&=& (2\pi)^{-n/2}t^{-\frac n{2\alpha}}\int_{\mr^n}\mx{e}^{i\frac{x}{t^{1/2\alpha}}\eta}
\mx{e}^{-|\eta|^{2\alpha}}\md\eta\\ \label{2.5}
&\triangleq& t^{-\frac n{2\alpha}}K\Big(\frac x{t^{1/2\alpha}}\Big).
\ee
Thus it is enough to consider the kernel function
\ben
K(x)=(2\pi)^{-n/2}\int_{\mr^n}\mx{e}^{ix\cdot\xi}\mx{e}^{-|\xi|^{2\alpha}}\md\xi.
\enn
It is obvious that $\mx{e}^{-|\xi|^{2\alpha}}\in L^1(\mr^n)$,
so
\be\label{2.8}
K(x)\in L^\infty(\mr^n)\cap C(\mr^n)
\ee
and by the Riemann-Lebesgue theorem,
$\lim_{|x|\rightarrow\infty}K(x)=0$.

Similarly, we have $|\xi|^\nu\mx{e}^{-|\xi|^{2\alpha}}\in L^1(\mr^n)$
and
\ben
(-\tr)^{\frac{\nu}2}K(x)\in L^\infty(\mr^n)\cap C_0(\mr^n)
\enn
for $\nu>0$, where $C_0(\mr^n)$ denotes the space of functions $f\in C(\mr^n)$
satisfying that $\lim\limits_{|x|\rightarrow\infty}f(x)=0$. In the same way, we have $\nabla K(x)\in L^\infty(\mr^n)\cap C_0(\mr^n)$
by the  fact $i\xi e^{-|\xi|^{\alpha}}\in \big( L^1(\Bbb R^n)\big)^n$.


\begin{lem}\label{lem2.1}
The kernel function $K(x)$ has the following point-wise estimate
\ben
|K(x)|\le C(1+|x|)^{-n-2\alpha},\;\; x\in \mr^n
\enn
for $\alpha>0$. Consequently one has
$$
K\in L^p(\mr^n),\;\;\ K_t\in L^p(\mr^n),\;\;0<t<\infty
$$
for any $1\le p\le\infty$.
\end{lem}

\begin{proof}
Define the invariant derivative operator
$$
L(x,D)=\frac{x\cdot\nabla_\xi}{i|x|^2}.
$$
Then we have
$$
L(x,D)\mx{e}^{ix\cdot\xi}=\mx{e}^{ix\cdot\xi}.
$$
The conjugate operator is
$$
L^*(x,D)=-\frac{x\cdot\nabla_\xi}{i|x|^2}.
$$
Thus we may write $K(x)$ as
\ben\no
K(x)&=&(2\pi)^{-n/2}\int_{\mr^n}\mx{e}^{ix\cdot\xi}L^*(\mx{e}^{-|\xi|^{2\alpha}})\md\xi\\ \no
&=&(2\pi)^{-n/2}\int_{\mr^n}\mx{e}^{ix\cdot\xi}
\rho(\frac\xi\delta)L^*(\mx{e}^{-|\xi|^{2\alpha}})\md\xi\\ \no
&&+(2\pi)^{-n/2}\int_{\mr^n}\mx{e}^{ix\cdot\xi}\Big(1-\rho(\frac\xi\delta)\Big)
L^*(\mx{e}^{-|\xi|^{2\alpha}})\md\xi\triangleq I+II,
\enn
where $\delta>0$ to be chosen later and  $\rho(\xi)$ is a $C^\infty_c(\mr^n)$-function satisfying
\ben
\rho(\xi)=\begin{cases}
1,\ &|\xi|\le 1;\\
0,\ &|\xi|>2.
\end{cases}
\enn
It is clear that
$$
|I|\le \frac{C}{|x|}\int_{|\xi|\le 2\delta}|\xi|^{2\alpha-1}\md\xi\le C|x|^{-1}\delta^{2\alpha+n-1}.
$$
To estimate $II$, take a sufficiently large natural number
$N>[2\alpha]+n$ and integrate by parts to obtain that
\ben
|II|&\le&(2\pi)^{-n/2}\int_{\mr^n}|\mx{e}^{ix\cdot\xi}(L^*)^{N-1}\Big((1-\rho(\frac\xi\delta))
L^*(\mx{e}^{-|\xi|^{2\alpha}})\Big)\md\xi\\ \nonumber
&\le& C|x|^{-N}\int_{|\xi|\ge\delta}\sum_{l=1}^N|\xi|^{2l\alpha-N}\mx{e}^{-|\xi|^{2\alpha}}\md\xi\\ \no
&&+C|x|^{-N}\sum_{k=1}^{N-1}C_k\delta^{-k}\int_{\delta\le |\xi|\le 2\delta}
\sum_{l=1}^{N-k}C_l|\xi|^{2l\alpha-(N-k)}\mx{e}^{-|\xi|^{2\alpha}}\md\xi\\ \no
&\le& C|x|^{-N}\int_{|\xi|\ge\delta}|\xi|^{2\alpha-N}\mx{e}^{-|\xi|^{2\alpha}}\md\xi
+C|x|^{-N}\int_{|\xi|\ge\delta}|\xi|^{2\alpha-N}|\xi|^{2\alpha(N-1)}\mx{e}^{-|\xi|^{2\alpha}}\md\xi\\
&&+C|x|^{-N}\sum_{k=1}^{N-1}\int_{\delta\le|\xi|\le 2\delta}\Big(|\xi|^{2\alpha-N}
\mx{e}^{-|\xi|^{2\alpha}}+|\xi|^{2\alpha(N-k)-N}\mx{e}^{-|\xi|^{2\alpha}}\Big)\md\xi.
\enn
In view of the facts that
$$\ds|\xi|^{2\alpha (N-1)}\mx{e}^{-|\xi|^{2\alpha}}\le C,\qquad
\ds|\xi|^{2\alpha (N-k-1)}\mx{e}^{-|\xi|^{2\alpha}}\le C
$$
for $k=1,\ 2,\ \cdots,\ N-1$, it is found that $|II|$ is dominated by
\ben
C|x|^{-N}\Big(\int_{|\xi|\ge\delta}|\xi|^{2\alpha-N}\md\xi+\int_{\delta\le|\xi|\le
2\delta}\delta^{2\alpha-N}\md\xi\Big)\le C|x|^{-N}\delta^{2\alpha-N+n}.
\enn
Thus it follows that
\ben
|K(x)|\le C|x|^{-1}\delta^{2\alpha+n-1}+C|x|^{-N}\delta^{2\alpha-N+n}.
\enn
Taking $\delta=|x|^{-1}$ gives
\ben
|K(x)|\le C|x|^{-n-2\alpha}.
\enn
This together with the boundedness of $K(x)$ (see (\ref{2.8})) completes the proof of the lemma.
\end{proof}

We now take the $\nu$-th derivative of the kernel $K(x)$ and have
$$
K^\nu(x)=(-\tr)^{\nu/2}K(x), \;\quad
K^\nu_t(x)=(-\tr)^{\nu/2}K_t(x).
$$
Then $K^\nu(x)$ can be split up into
\ben
K^\nu(x)&=&(2\pi)^{-n/2}\int_{\mr^n}\mx{e}^{ix\cdot\xi}\rho(\xi/\delta)
|\xi|^\nu\mx{e}^{-|\xi|^{2\alpha}}\md\xi\\
&&\quad+(2\pi)^{-n/2}\int_{\mr^n}\mx{e}^{ix\cdot\xi}(1-\rho(\xi/\delta))
|\xi|^\nu\mx{e}^{-|\xi|^{2\alpha}}\md\xi\\
&\triangleq& I+II.
\enn
Clearly,
\ben
|I|\le C\int_{|\xi|\le 2\delta}\delta^\nu\md\xi\le C\delta^{n+\nu}.
\enn
To estimate $II$ we use the technique of invariant derivatives together with
integration by parts to obtain that
\ben
|II|\le C\int_{\mr^n}\Big|\mx{e}^{ix\cdot\xi}(L^*)^N\big((1-\rho(\xi/\delta))
|\xi|^\nu\mx{e}^{-|\xi|^{2\alpha}}\big)\Big|\md\xi.
\enn
Arguing similarly as the proof of Lemma \ref{lem2.1} we have
\ben\nonumber
|II|&\le& C|x|^{-N}\Big(\int_{|\xi|\ge\delta}|\xi|^{\nu-N}\md\xi
+\int_{\delta\le |\xi|\le 2\delta}\delta^{\nu-N}\md\xi\Big)\\
&\le& C|x|^{-N}\delta^{\nu-N+n}.
\enn
Taking $\delta=|x|^{-1}$ leasd to the estimate
$$
|K^\nu(x)|\le C|x|^{-\nu-n}.
$$
Thus we have the following lemma.

\begin{lem}\label{lem2.2}
The kernel function $K^\nu(x)$ has the following pointwise estimate
\ben
|K^\nu(x)|\le C(1+|x|)^{-n-\nu},\;\; x\in\mr^n
\enn
for $\nu>0$. Consequently one has
$$
K^\nu\in L^p(\mr^n),\;\;\; K^\nu_t\in L^p(\mr^n),\;\; 0<t<\infty
$$
for any $1\le p\le\infty$.
\end{lem}

\begin{rem}\label{r2.1} {\rm
(i)  Thank to  $i\xi e^{-|\xi|^{\alpha}}\in \big( L^1(\Bbb R^n)\big)^n$, one has  by the same argument of  Lemma \ref{lem2.2} that
$$
|\nabla K(x)|\le C(1+|x|)^{-n-1},
$$
and
$$
\nabla K(x),\;\; \nabla K_t(x)\in L^p(\mr^n),\ 0<t<\infty
$$
for any $1\le p\le \infty$.

(ii) Similar  to (\ref{2.5}) the kernel function $K^\nu_t(x)$ satisfies the same scaling as follows:
\be\label{2.24}
K^\nu_t(x)=t^{-\frac\nu{2\alpha}}t^{-\frac n{2\alpha}}K^\nu\Big(\frac x{t^{1/2\alpha}}\Big).
\ee
}
\end{rem}

In Proposition \ref{prop2.1} below we give another characterization of the negative index
homogeneous Besov space $\dot{B}^s_{p,q}(\mr^n),$ employing the pointwise estimate
of the kernel $K(x)$ in Lemma \ref{lem2.1} and the boundedness of the fractional
power dissipative operator semigroup $S_\alpha(t)$ on the space $L^p(\mr^n),$ where $s<0$.
The idea essentially comes from \cite{Peetre} (see also \cite{MYang}).
For the case $n=2$ and $1/2<\alpha\le 1$ the reader is also referred to \cite{Zhang}.
For completeness we give a proof of Proposition \ref{prop2.1} here for any $n\in\mb{N}$
and $0<\alpha<\infty$. We first recall the definition of homogeneous Besov spaces.

Choose a radial bump function $\widehat{\psi}(\xi)\in C^\infty_c(\mathbb{R}^n)$ such that
\ben\nonumber
\widehat{\psi}(\xi)=\left\{\begin{array}{ll}
  1,\quad \mbox{if}\ |\xi|\le 1,\\[0.2cm]
  \mbox{smooth},\quad \mbox{if}\ 1<|\xi|<2,\\[0.2cm]
  0,\quad \mbox{if}\ |\xi|\ge 2,
 \end{array}\right.
\mbox{and}\quad  0\le \widehat{\psi}(\xi)\le 1,
\enn
where $\widehat{\psi}(\xi)$ denotes the Fourier transform of $\psi(x)$.
Set $\widehat{\phi}(\xi)=\widehat{\psi}(\xi)-\widehat{\psi}(2\xi)$ and let
$\widehat{\phi}_j(\xi)=\widehat{\phi}(2^{-j}\xi),\ \xi\neq 0$ for $j\in\mathbb{Z},$
$\widehat{\psi}_j(\xi)=\widehat{\psi}(2^{-j}\xi)$ for $j\in\mathbb{Z}$.
Let $\tr_jf=\phi_j*f,$ $S_jf=\psi_j*f$. Then for any $f\in L^2(\mathbb{R}^n)$
we have the following Littlewood-Paley decomposition
\ben
f(x)=\sum^\infty_{j=-\infty}\phi_j*f(x),
\enn
where the sum is taken in the $L^2(\mathbb{R}^n)$ sense.

The homogeneous Besov space $\dot{B}^s_{p,q}$ is defined by the
dyadic decomposition as
$$
\dot{B}^s_{p,q}=\{f\in\mathcal{Z}'(\mathbb{R}^n)\,|\,\|f\|_{\dot{B}^s_{p,q}}<\infty\},
$$
where
$$
\|f\|_{\dot{B}^s_{p,q}}=
\begin{cases}\ds\left(\sum^\infty_{j=-\infty}2^{jsq}\|\phi_j*f\|^q_p\right)^{1/q},&1\le q<\infty\\
             \ds\sum_j2^{js}\|\phi_j*f\|_p,&q=\infty
\end{cases}
$$
is the norm of $\dot{B}^s_{p,q}$ and $\mathcal{Z}'(\mathbb{R}^n)$ denotes the dual space of
$$
\mathcal{Z}(\mathbb{R}^n)=\{f\in\mathcal{S}(\mathbb{R}^n)\,\bigg|D^\al\hat{f}(0)=0,\;
\mbox{for any multi-index}\;\al\in\mathbb{N}^n\}
$$
and can be identified by the quotient space $\mathcal{S}'/\mathcal{P}$ with
the polynomial $\mathcal{P}$. See \cite{B-L}, \cite{M1} and \cite{Tr} for details.

\begin{prop}\label{prop2.1}
Let $1\le p,\ q\le \infty$, $s<0$ and assume that $n\in\mb{N}$ and $0<\alpha<\infty.$
Then $\ds f\in\dot{B}^s_{p,q}(\mr^n)$ if and only if
\be\label{2.27}
\begin{cases}\ds\bigg(\int_0^\infty\Big(t^{-\frac s{2\al}}\|S_\al(t)f\|_p\Big)^q
                \frac{\md t}t\bigg)^{1/q}<\infty,&1\le q<\infty,\\
                \ds\sup_{t>0}t^{-\frac s{2\al}}\|S_\al(t)f\|_p,&q=\infty.
\end{cases}
\ee
\end{prop}

\begin{proof} We only consider the case $1\le q<\infty.$ The case $q=\infty$ can be shown similarly.
We first prove that
\ben
\bigg(\sum^\infty_{j=-\infty}2^{jsq}\|\tr_jf\|^q_p\bigg)^{1/q}
\le C\bigg(\int_0^\infty\Big(t^{-\frac s{2\alpha}}\|S_\alpha(t)f\|_p\Big)^q\frac{\md t}t\bigg)^{1/q}.
\enn
In fact, let
\ben
\Phi_j(x)&=&\mc{F}^{-1}\Big(\widehat{\phi}(\frac\xi{2^j})\mx{e}^{(2^{-j}|\xi|)^{2\alpha}}\Big)(x),\\
h_t(x)&=&\mc{F}^{-1}\Big(\mx{e}^{-(t|\xi|)^{2\alpha}}\Big)(x).
\enn
Then by the definition of $\tr_j$ one has $\tr_jf=\Phi_j*h_{2^{-j}}*f(x)$.

By the Young inequality we get
\ben
\|\tr_jf\|p\le\|\Phi_j\|_1\|h_{2^{-j}}*f\|_p\le C\|h_{2^{-j}}*f\|_p,
\enn
where we have used the fact that
\ben
\|\Phi_j(x)\|_1=\int_{\mr^n}|\mc{F}^{-1}\Big(\mx{e}^{|\xi|^{2\alpha}}
\widehat{\phi}(\xi)\Big)(x)|\md x<\infty.
\enn
A direct calculation shows that
\ben\nonumber
h_{2^{-j}}*f(x)&=&(2\pi)^{-n/2}\int_{\mr^n}\mx{e}^{-(2^{-j}|\xi|)^{2\alpha}}\widehat{f}(\xi)
\mx{e}^{ix\cdot\xi}\md\xi\\ \nonumber
&=&(2\pi)^{-n/2}\int_{\mr^n}\mx{e}^{-|\xi|^{2\alpha}(2^{-2\alpha j}-t^{2\alpha})}
\mx{e}^{-(t|\xi|)^{2\alpha}}\widehat{f}(\xi)\mx{e}^{ix\cdot\xi}\md\xi\\
&=& S_\alpha(2^{-2\alpha j}-t^{2\alpha})(h_t*f)(x).
\enn
Thus it follows that
\be\label{2.34}
\|h_{2^{-j}}*f\|_p\le C\|h_t*f\|_p
\ee
for any $t\in [2^{-j-1},2^{-j}]$, which implies that
\ben\nonumber
\sum_{j=-\infty}^\infty2^{sjq}\|\tr_jf\|_q^q&\le& C\sum_{j=-\infty}^\infty
\int^{2^{-j}}_{2^{-j-1}}\Big(t^{-s}\|h_{2^{-j}}*f\|_p\Big)^q\frac{\md t}t\\\nonumber
&\le& C\int^\infty_0\Big(t^{-s}\|h_t*f\|_p\Big)^q\frac{\md t}t\\
&\le& C\int^\infty_0\Big(t^{-\frac s{2\alpha}}\|S_\alpha(t)*f\|_p\Big)^q\frac{\md t}t,
\enn
where we have used  the fact that $h_t*f(x)=S_\alpha(t^{2\alpha})f(x)$.

We now prove that
\be\label{2.36}
\bigg(\int_0^\infty\Big(t^{-\frac s{2\alpha}}\|S_\alpha(t)f\|_p\Big)^q\frac{\md t}t\bigg)^{1/q}
\le C\bigg(\sum^\infty_{j=-\infty}2^{jsq}\|\tr_jf\|^q_p\bigg)^{1/q}.
\ee
In fact, for any $j\in\mb{Z}$ one has the decomposition
\ben
h_{2^{-j}}*f(x)=\sum_{k=-\infty}^\infty h_{2^{-j}}*\tr_{k+j}f(x).
\enn
Arguing similarly as in deriving (\ref{2.34}) one has
\ben
\|h_t*f(x)\|_p\le C\|h_{2^{-j}*f(x)}\|_p
\enn
for any $t\in [2^{-j},2^{-j+1}].$
The left-hand side of the estimate (\ref{2.36}) can be estimated as follows:
\be\nonumber
\int^\infty_0\Big(t^{-\frac s{2\alpha}}\|S_\alpha(t)f\|_p\Big)^q\frac{\md t}t
&=&2\alpha\int^\infty_0\Big(t^{-s}\|h_t*f\|_p\Big)^q\frac{\md t}t\\ \nonumber
&\le& C\sum_{j=-\infty}^\infty\int^{2^{-j+1}}_{2^{-j}}\Big(2^{js}
\|h_{2^{-j}}*f\|_p\Big)^q\frac{\md t}t\\ \label{2.39}
&\le& C\sum_{j=-\infty}^\infty\Big(2^{js}\sum_{k=-\infty}^\infty\|h_{2^{-j}}*\tr_{k+j}f\|_p\Big)^q.
\ee
If we can show that
\be\label{2.41}
\|h_{2^{-j}}*\tr_{k+j}f\|_p\le 2^{ks}\|\tr_{k+j}f\|_p
\ee
for any $s<0,$ then taking $s_1<s<s_0<0$ we have by using
the Minkowski inequality that the right-hand side of (\ref{2.39}) is bounded above by
\ben\nonumber
&& C\sum_{j=-\infty}^\infty\Big(2^{js}\sum_{k=-\infty}^02^{ks_0}\|\tr_{k+j}f\|_p\Big)^q
   +C\sum_{j=-\infty}^\infty\Big(2^{js}\sum_{k=1}^\infty2^{ks_1}\|\tr_{k+j}f\|_p\Big)^q\\ \no
&&\quad\le C\bigg(\sum^0_{-\infty}2^{k(s_0-s)}\Big(\sum_{j=-\infty}^\infty2^{(k+j)sq}
  \|\tr_{k+j}f\|^q_p\Big)^{1/q}\bigg)^q\\ \nonumber
&&\qquad+ C\bigg(\sum^\infty_{k=1}2^{-k(s-s_1)}\Big(\sum_{j=-\infty}^\infty2^{(k+j)sq}
  \|\tr_{k+j}f\|^q_p\Big)^{1/q}\bigg)^q\\
&&\quad\le C\sum_{j=-\infty}^\infty2^{sjq}\|\tr_jf\|_p^q,
\enn
which completes the proof of Proposition \ref{prop2.1}. We now prove the estimate (\ref{2.41}).
Note first that
\ben\nonumber
h_{2^{-j}}*\tr_{k+j}f(x)&=&\frac1{(2\pi)^{n/2}}\int_{\mr^n}\mx{e}^{-(2^{-j}|\xi|)^{2\alpha}}
\widehat\phi(2^{-k-j}\xi)\widehat{f}(\xi)\mx{e}^{ix\cdot\xi}\md \xi\\ \nonumber
&=&\frac{2^{ks}}{(2\pi)^{n/2}}\int_{\mr^n}\mx{e}^{-(2^{-j}|\xi|)^{2\alpha}}\Big(\frac{2^j}{|\xi|}\Big)^{s}
     \Big(\frac{|\xi|}{2^{k+j}}\Big)^{s}\widetilde{\phi}(2^{-k-j}\xi)
     \widehat\phi(2^{-k-j}\xi)\widehat{f}(\xi)\mx{e}^{ix\cdot\xi}\md\xi,
\enn
where $\widetilde{\phi}(2^{-j}\xi)=\widehat\phi(2^{-j+1}\xi)+\widehat\phi(2^{-j}\xi)
+\widehat\phi(2^{-j-1}\xi).$ Since
\ben\nonumber
\Big\|\mc{F}^{-1}\Big(\mx{e}^{-(2^{-j}|\xi|)^{2\alpha}}\Big(\frac{2^j}{|\xi|}\Big)^{s'}\Big)\Big\|_{L^1}
&=&\Big\|2^{jn}\mc{F}^{-1}(\mx{e}^{-|\xi|^{2\alpha}|\xi|^{-s}})(2^jx)\Big\|_{L^1}\\
&=&\Big\|\mc{F}^{-1}(\mx{e}^{-|\xi|^{2\alpha}|\xi|^{-s}})\Big\|_{L^1}<\infty,\\
\Big\|\mc{F}^{-1}\Big(\left(\frac{|\xi|}{2^{k+j}}\right)^{s}
\widetilde{\phi}(\frac{\xi}{2^{k+j}})\Big)\Big\|_{L^1}&=&
\|\mc{F}^{-1}(|\xi|^{s}\widetilde{\phi}(\xi))\|_{L^1}<\infty,
\enn
the estimate (\ref{2.41}) follows easily from the Young inequality.
\end{proof}


\section{Space-time estimates for the linear equation}
\setcounter{equation}{0}

In this section we discuss the space-time estimates of solutions to the Cauchy problem of the linear
fractional power dissipative equation
\begin{eqnarray}\label{3.1}
\begin{cases}
u_t+(-\tr)^\alpha u=f(t,x),\ &(t,x)\in [0,\infty)\times\mr^n;\\
u(0)=\varphi(x),\ & x\in \mr^n.
\end{cases}
\end{eqnarray}
By Duhamel's principle, the solution to the problem (\ref{3.1}) can be written in the
integral form as
\be\label{3.3}
u(x,t)=S_\alpha(t)\varphi(x)+\int^t_0S_\alpha(t-\tau)f(\tau,x)\md\tau
\triangleq S_\alpha(t)\varphi(x)+(\mb{G}f)(t,x).
\ee

%
We first consider the space-time estimates for the homogeneous part of the solution $u$
given in the integral form \ref{3.3}.

\begin{lem}\label{lem3.1}
Let $1\le r\le p\le\infty$ and let $\varphi\in L^r(\mr^n)$. Then the
homogeneous part of the solution (\ref{3.3}) satisfies the estimates
\be\label{3.5}
\|S_\alpha(t)\varphi(x)\|_p&\le& Ct^{-\frac n{2\alpha}(\frac1r-\frac1p)}\|\varphi\|_{L^r},\\ \label{3.6}
\|(-\tr)^{\nu/2}S_\alpha(t)\varphi(x)\|_p
&\le& Ct^{-\frac\nu{2\alpha}-\frac n{2\alpha}(\frac1r-\frac1p)}\|\varphi\|_{L^r}
\ee
for $\alpha>0$ and $\nu>0$.
\end{lem}

\begin{proof} It follows from the Young inequality combined with scaling property of the kernel $K_t$.
\end{proof}

To derive the space-time estimates of the homogeneous part of the solution $u$
given in \ref{3.3}, we need to introduce the following definition on
admissible triplets and generalized admissible triplets for the
fractional power dissipative equation. For the corresponding definition
for parabolic equations the reader is referred to \cite{M-G,Miao1,M-Z}.

\begin{defin}\label{def3.1}
The triplet $(q,p,r)$ is called an admissible triplet (for the
fractional power dissipative equation) if
\ben
\frac1q=\frac n{2\alpha}\bigg(\frac1r-\frac1p\bigg),
\enn
where
\ben
1<r\le p<\begin{cases}
\frac{nr}{n-2\alpha},\ &\mx{for } n>2\alpha,\\
\infty,\ &\mx{for } n\le 2\alpha.
          \end{cases}
\enn
\end{defin}

\begin{defin}\label{def3.2}
The triplet $(q,p,r)$ is called a generalized admissible triplet (for the
fractional power dissipative equation) if
\ben
\frac1q=\frac n{2\alpha}\bigg(\frac1r-\frac1p\bigg),
\enn
where
\ben
1<r\le p<\begin{cases}
\frac{nr}{n-2\alpha r},\ &\mx{for } n>2r\alpha,\\
\infty,\ &\mx{for } n\le 2r\alpha.
          \end{cases}
\enn
\end{defin}

Let B be a Banach space and let $I=[0,T)$. We define the time-weighted
space-time Banach space $\mc{C}_\sigma(I;B)$ and the corresponding homogeneous space
$\dot{\mc{C}}_\sigma(I;B)$ as follows
\ben
\mc{C}_\sigma(I;B)&=&\{f\in C(I;B)\,\big|\,\|f;\mc{C}_\sigma(I;B)\|
  =\sup_{t\in I}t^{\frac1\sigma}\|f\|_B<\infty\},\\
\dot{\mc{C}}_\sigma(I;B)&=&\{f\in \mc{C}_\sigma(I;B)\,\big|\,
\lim_{t\rightarrow 0^+}t^{\frac1\sigma}\|f\|_B=0\}. \enn In this
paper the Banach space $B$ is taken to be $L^p(\mr^n)$ with $1< p<
\infty.$

With the above definitions we now have the following results on the space-time estimates
for the homogeneous part of the solution $u$ given in (\ref{3.3}).
These estimates can be proved by following \cite{Giga} (see also \cite{M-G}).
Here we give a proof for completeness.

\begin{lem}\label{lem3.2}
(i) Let $(q,p,r)$ be any admissible triplet and let $\varphi\in L^r(\mr^n)$.
Then $\ds S_\alpha(t)\varphi\in L^q(I;L^p(\mr^n))
\cap C_b(I;L^r(\mr^n))$ with the estimate
\be\label{3.13}
\|S_\alpha(t)\varphi(x)\|_{L^q(I;L^p)}\le C\|\varphi\|_{L^r},
\ee
for $0<T\le\infty$, where $C$ is a positive constant.

(ii) Let $(q,p,r)$ be any generalized admissible triplet. For any $\varphi\in L^r(\mr^n)$
we have $\ds S_\alpha(t)\varphi\in\mc{C}_q(I;L^p(\mr^n))\cap C_b(I;L^r(\mr^n))$ and
\be\label{3.15}
\|S_\alpha(t)\varphi\|_{\mc{C}_q(I;L^p)}\le C\|\varphi\|_{L^r}.
\ee
Hereafter, for a Banach space $X$ we denote by $C_b(I;X)$ the space of bounded continuous functions
from $I$ to $X$.
\end{lem}

\begin{proof}
The statement (ii) follows easily from Lemma \ref{lem3.1}. So we only need to prove (i).
For the case $p=r,\, q=\infty$, the estimate (\ref{3.13}) is true from Lemma \ref{lem3.1}.
We now consider the case $p\not=r$. Let
\ben
U(t)\varphi=\|S_\alpha(t)\varphi\|_p.
\enn
Then, and since $(q,p,r)$ is an admissible triplet, we deduce by Young's inequality that
\ben
U(t)\varphi\le Ct^{-\frac1q}\|\varphi\|_{L^r}.
\enn
It is easy to see that
\ben\nonumber
\mu\{t:\;|U(t)\varphi|>\tau\}&\le& \mu\{t:\;Ct^{-\frac1q}\|\varphi\|_{L^r}>\tau\}
=\mu\bigg\{t:\;t<\Big(\frac{C\|\varphi\|_{L^r}}\tau\Big)^q\bigg\}\\
&\le&\Big(\frac{C\|\varphi\|_{L^r}}\tau\Big)^q,
\enn
which implies that $U(t)$ is a weak type $(r,q)$ operator.

On the other hand, by Lemma \ref{lem3.1} $U(t)$ is sub-additive and satisfies that
\ben
U(t)\varphi=\|S_\alpha(t)\varphi\|_p\le C\|\varphi\|_p
\enn
for $r\le p\le\infty$, which means that $U(t)$ is a $(p,\infty)$ operator.
Since for any admissible triplet $(p,q,r)$ we can always find another admissible triplet
$(p,q_1,r_1)$ such that
$$
q_1< q<\infty,\quad\ r_1<r<p
$$
and
$$
\frac1q=\frac\theta{q_1}+\frac{1-\theta}\infty,\quad
\frac1r=\frac\theta{r_1}+\frac{1-\theta}p,
$$
then the Marcinkiewicz interpolation theorem (see \cite{M1} or \cite{S-W}) implies
that $U(t)$ is a strong $(r,q)$-type operator. The estimate (\ref{3.13}) thus follows,
and the proof of Lemma \ref{lem3.2} is complete.
\end{proof}

We now derive the space-time estimates of the non-homogeneous part $\mb{G}f$ of the solution $u$
given in (\ref{3.3}). For the case when $\alpha$ is a positive integer, see also \cite{M1,M-Z}.

\begin{lem}\label{lem3.3}
For $b>0$ and $T>0$ let $r_0={nb}/({2\alpha}),$ $I=[0,T)$.
Assume that $r\ge r_0>1$ and that $(q,p,r)$ is an admissible triplet satisfying that $p>b+1$.

(i) If $\ds f\in L^\frac q{b+1}(I;L^{\frac p{b+1}}(\mr^n))$, then
\ben
\|\mb{G}f\|_{L^\infty(I;L^r)}\le CT^{1-\frac{nb}{2r\alpha}}
\|f\|_{L^{\frac q{b+1}}(I;L^{\frac p{b+1}})}
\enn
for $p<r(1+b)$, and
\ben
\|\mb{G}f\|_{L^\infty(I;L^r)}\le CT^{1-\frac{nb}{2r\alpha}}
\||f|^{\frac1{b+1}}\|^{\theta(b+1)}_{L^\infty(I;L^r)}
\||f|^{\frac1{b+1}}\|^{(1-\theta)(b+1)}_{L^q(I;L^p)}
\enn
for $p\ge r(b+1)$, where $\ds\theta=\frac{p-r(b+1)}{(b+1)(p-r)}$.

(ii) If $\ds f\in L^\frac q{b+1}(I;L^{\frac p{b+1}}(\mr^n))$, then
\ben
\|\mb{G}f\|_{L^q(I;L^p)}\le CT^{1-\frac{nb}{2r\alpha}}
\|f\|_{L^{\frac q{b+1}}(I;L^{\frac p{b+1}})}
\enn
for $p<r(b+1),$ and
\ben
\|\mb{G}f\|_{L^q(I;L^p)}\le CT^{1-\frac{nb}{2r\alpha}}
\||f|^{\frac1{b+1}}\|^{\theta(b+1)}_{L^\infty(I;L^r)}
\||f|^{\frac1{b+1}}\|^{(1-\theta)(b+1)}_{L^q(I;L^p)}
\enn
for $p\ge r(b+1),$ where $\theta$ is the same as in (i).
\end{lem}

\begin{proof}
We first prove (i). Consider first the case when $p<r(b+1)$. Using Young's inequality one has
\ben\nonumber
\|\mb{G}f\|_{L^\infty(I;L^r)}&\le&C\int^t_0(t-s)^{-\frac n{2\alpha}(\frac{b+1}p-\frac1r)}
   \|f(s,x)\|_{L^{\frac p{b+1}}}\md s\\ \nonumber
&\le& C\bigg(\int^t_0(t-s)^{-\frac n{2\alpha}(\frac{b+1}p-\frac1r)\chi}\md s\bigg)^{\frac1\chi}
   \|f\|_{L^{\frac q{b+1}}(I;L^{\frac p{b+1}})}\\
&\le& CT^{1-\frac{nb}{2\alpha r}}\|f\|_{L^{\frac q{b+1}}(I;L^{\frac p{b+1}})},
\enn
where $\ds\frac1{\chi}=1-\frac{b+1}q$ and $C=C(n,p,r,b)$ depends only on $n,\ p,\ r,\ b$.

For the case $p\ge r(b+1)$, by means of the Riesz interpolation theorem (see e.g. \cite{S-W,M1} or
\cite{B-L}) and the H\"older inequality we have
\ben\nonumber
\|\mb{G}f\|_{L^\infty(I;L^r)} &\le& \int^t_0\||f(s,x)|^{\frac1{b+1}}\|^{b+1}_{r(b+1)}\md s\\ \no
&=& C\int^t_0\||f(s,x)|^{\frac1{b+1}}\|^{(b+1)\theta}_{L^r}
\||f(s,x)|^{\frac1{b+1}}\|^{(b+1)(1-\theta)}_p\md s\\ \nonumber
&\le& CT^{1-\frac{(b+1)(1-\theta)}q}\||f|^{\frac1{b+1}}\|^{(b+1)\theta}_{C(I;L^r)}
\||f|^{\frac1{b+1}}\|^{(b+1)(1-\theta)}_{L^q(I;L^p)}\\
&\le& CT^{1-\frac{nb}{2\alpha r}}\||f|^{\frac1{b+1}}\|^{(b+1)\theta}_{C(I;L^r)}
\||f|^{\frac1{b+1}}\|^{(b+1)(1-\theta)}_{L^q(I;L^p)},
\enn
where $\theta$ satisfies $\ds\frac1{r(b+1)}=\frac\theta r+\frac{1-\theta}p,$
the index of the H\"older inequality is $\ds 1=\frac{(1+b)(1-\theta)}q+\frac1\chi$,
and use has been made of the fact that
\be\no
1-\frac{(b+1)(1-\theta)}q&<&1-\frac{b+1}q+\frac{n(b+1)}{2\alpha}\bigg(\frac1{r(b+1)}-\frac1p\bigg)\\
 \label{Tindex}
&=& 1-\frac{nb}{2\alpha r}.
\ee
We now prove (ii). For the case $p<r(b+1)$ we have by Young's inequality that
\ben\nonumber
\|\mb{G}f\|_{L^q(I;L^p)} &\le& C\left\|\int^t_0(t-s)^{-\frac n{2\alpha}
(\frac {b+1}p-\frac1p)}\|f(s,x)\|_{L^{\frac p{b+1}}}\md s\right\|_{L^q}\\ \nonumber
&\le& C\bigg(\int^T_0t^{-\frac{nb}{2{\alpha p}}\chi}\bigg)^{\frac1\chi}
\|f\|_{L^{\frac q{b+1}}(I;L^{\frac p{b+1}})}\\
&\le& CT^{1-\frac{nb}{2\alpha r}}\|f\|_{L^{\frac q{b+1}}(I;L^{\frac p{b+1}})},
\enn
where $\ds 1+\frac1q=\frac{1+b}q+\frac1\chi$.
For $p\ge r(b+1)$, arguing similarly as in the proof of (i) gives
\ben\nonumber
\|\mb{G}f\|_{L^q(I;L^p)}&\le& C\left\|\int^t_0(t-s)^{-\frac n{2\alpha}(\frac1r-\frac1p)}
\||f(s,x)|^{\frac1{b+1}}\|^{b+1}_{r(b+1)}\md s\right\|_q\\ \nonumber
&\le& C\left\|\int^t_0(t-s)^{-\frac n{2\alpha}(\frac1r-\frac1p)}
\||f(s,x)|^{\frac1{b+1}}\|^{(b+1)\theta}_{L^r}
\||f(s,x)|^{\frac1{b+1}}\|^{(b+1)(1-\theta)}_p\md s\right\|_q\\ \nonumber
&\le& C\bigg(\int^T_0t^{-\frac n{2\alpha}(\frac1r-\frac1p)\chi}\md t\bigg)^{\frac1\chi}
\||f|^{\frac1{b+1}}\|^{\theta(b+1)}_{C(I;L^r)}
\||f|^{\frac1{b+1}}\|^{(b+1)(1-\theta)}_{L^q(I;L^p)}\\
&\le& CT^{1-\frac{nb}{2\alpha r}}\||f|^{\frac1{b+1}}\|^{\theta(b+1)}_{C(I;L^r)}
\||f|^{\frac1{b+1}}\|^{(b+1)(1-\theta)}_{L^q(I;L^p)},
\enn
where $\theta$ and $\chi$ satisfy that
$$
\frac1{r(b+1)}=\frac\theta r+\frac{1-\theta}p,\qquad
1+\frac1q=\frac{(b+1)(1-\theta)}q+\frac1\chi,
$$
which is meaningful by the fact that $r<r(b+1)<p.$
\end{proof}

Arguing similarly in the proof of Lemma \ref{lem3.3} we can derive the estimates in
the spaces $\ds\mc{C}_q(I;L^p(\mr^n))$ and $C_b(I;L^r(\mr^n))$ of the non-homogeneous term.
In fact, for the case $p<r(b+1)$ (which implies $q>b+1$), one has by Lemma \ref{lem3.1} that
\ben\nonumber
\|\mb{G}f\|_{L^\infty(I;L^r)} &\le& C\int^t_0(t-s)^{-\frac n{2\alpha}(\frac{b+1}p-\frac1r)}
   \|f(s,x)\|_{L^{\frac p{b+1}}}\md s\\ \nonumber
&\le& C\int^t_0(t-s)^{-\frac n{2\alpha}(\frac{b+1}p-\frac1r)}s^{-\frac{b+1}q}\md s
\|f\|_{\mc{C}_{\frac q{b+1}}(I;L^{\frac p{b+1}})}\\
&\le& CT^{1-\frac{nb}{2r\alpha}}\|f\|_{\mc{C}_{\frac q{b+1}}(I;L^{\frac p{b+1}})},
\enn
where $C=C(n,p,r,b)$ depends only on $n,\ p,\ r,\ b$.
Making use of the space-time estimates for the heat equation (cf. \cite{M-G}) and
Young's inequality we get
\ben\nonumber
\|\mb{G}f\|_{\mc{C}_q(I;L^p)} &\le& C\sup_{t\in I}t^{\frac1q}
\int^t_0(t-s)^{-\frac n{2\alpha}(\frac {b+1}p-\frac1p)}\|f(s,x)\|_{\frac p{b+1}}\md s\\ \no
 &\le& C\sup_{t\in I }t^{\frac1q}\int^t_0(t-s)^{-\frac{nb}{2p\alpha}}s^{-\frac{b+1}q}\md s
 \|f\|_{\mc{C}_{\frac q{b+1}}(I;L^{\frac p{b+1}})}\\
&\le& CT^{1-\frac{nb}{2r\alpha}}\|f\|_{\mc{C}_{\frac q{b+1}}(I;L^{\frac p{b+1}})}.
\enn

For the case $p\ge r(b+1)$, we use the Riesz interpolation theorem (see \cite{S-W,M1} or \cite{B-L})
to get, on noting the definition of the space $\mc{C}_q(I;L^p(\mr^n)),$ that for any $0<t\le T$
\be\no
\|\mb{G}f\|_{L^\infty(I;L^r)} &\le& \int^t_0\||f(s,x)|^{\frac1{b+1}}\|^{b+1}_{r(b+1)}\md s\\ \no
&=& C\int^t_0\||f(s,x)|^{\frac1{b+1}}\|^{(b+1)\theta}_{L^r}\||f(s,x)|^{\frac1{b+1}})\|_p^{(b+1)
    (1-\theta})\md s\\ \no
&\le&C\||f|^{\frac1{b+1}}\|^{(b+1)\theta}_{C(I;L^r)}
\||f|^{\frac1{b+1}}\|^{(1+b)(1-\theta)}_{\mc{C}_q(I;L^p)}
\int^t_0s^{-\frac1{q}(b+1)(1-\theta)}\md s\\ \nonumber
&\le& C\||f|^{\frac1{b+1}}\|^{(b+1)\theta}_{C(I;L^r)}
\||f|^{\frac1{b+1}}\|^{(1+b)(1-\theta)}_{\mc{C}_q(I;L^p)}T^{1-\frac1q(b+1)(1-\theta)}\\ \label{West1}
&\le& C T^{1-\frac{nb}{2\alpha r}}\||f|^{\frac1{b+1}}\|^{(b+1)\theta}_{C(I;L^r)}
\||f|^{\frac1{b+1}}\|^{(1+b)(1-\theta)}_{\mc{C}_q(I;L^p)},
\ee
where $\theta$ satisfies that $\ds\frac1{r(b+1)}=\frac\theta r+\frac{1-\theta}p$.
To get the estimate in the time weighted space $\ds\mc{C}_q(I;L^p(\mr^n))$ we make use of the Riesz
interpolation theorem again and obtain that
\be\nonumber
\|\mb{G}f\|_{\mc{C}_{q}(I;L^p)}&\le& C\sup_{t\in I}t^{\frac1q}
\int^t_0(t-s)^{-\frac n{2\alpha}(\frac1r-\frac1p)}\||f(s,x)|^{\frac1{b+1}}\|^{b+1}_{r(b+1)}\md s\\ \no
&\le& C\sup_{t\in I}t^{\frac1q}\int^t_0(t-s)^{-\frac n{2\alpha}(\frac1r-\frac1p)}
\||f(s,x)|^{\frac1{b+1}}\|^{(b+1)\theta}_{L^r}\||f(s,x)|^{\frac1{b+1}}\|^{(b+1)(1-\theta)}_p\md s\\ \no
&\le& C\sup_{t\in I}t^{\frac1q}\int^t_0(t-s)^{-\frac n{2\alpha}(\frac1r-\frac1p)}
s^{-\frac{(b+1)(1-\theta)}q}\md s\||f|^{\frac1{b+1}}\|^{(b+1)\theta}_{C(I;L^r)}
\||f|^{\frac1{b+1}}\|^{(b+1)(1-\theta)}_{\mc{C}_q(I;L^p)}\\ \label{West2}
&\le& CT^{1-\frac{nb}{2r\alpha}}\||f|^{\frac1{b+1}}\|^{(b+1)\theta}_{C(I;L^r)}
\||f|^{\frac1{b+1}}\|^{(b+1)(1-\theta)}_{\mc{C}_q(I;L^p)},
\ee
where $\theta$ is the same as in (\ref{West1}).
Thus we have obtained the following results.

\begin{lem}\label{lem3.4}
For $b>0$ and $T>0$, let $r_0={nb}/({2\alpha}),$ $I=[0,T)$. Assume that $r\ge r_0>1$.
Let $(q,p,r)$ be any generalized admissible triplet satisfying that $p>b+1$.

(i) If $\ds f\in\mc{C}_{\frac q{b+1}}(I;L^{\frac p{b+1}}(\mr^n))$, then
\ben
\|\mb{G}f\|_{L^\infty(I;L^r)}\le CT^{1-\frac{nb}{2r\alpha}}
\|f\|_{\mc{C}_{\frac q{b+1}}(I;L^{\frac p{b+1}})},
\enn
for $p< r(1+b),$ and
\ben
\|\mb{G}f\|_{L^\infty(I;L^r)}\le CT^{1-\frac{nb}{2r\alpha}}
\||f|^{\frac1{b+1}}\|^{\theta(b+1)}_{L^\infty(I;L^r)}
\||f|^{\frac1{b+1}}\|^{(1-\theta)(b+1)}_{\mc{C}_{q}(I;L^p)}
\enn
for $p\ge r(b+1)$, where $\ds\theta=\frac{p-r(b+1)}{(b+1)(p-r)}$.

(ii) If $\ds f\in \mc{C}_{\frac q{b+1}}(I;L^{\frac p{b+1}}(\mr^n))$, then
\ben
\|\mb{G}f\|_{\mc{C}_q(I;L^p)}\le CT^{1-\frac{nb}{2r\alpha}}
\|f\|_{\mc{C}_{\frac q{b+1}}(I;L^{\frac p{b+1}})}
\enn
for $p<r(b+1)$, and
\ben
\|\mb{G}f\|_{\mc{C}_q(I;L^p)}\le CT^{1-\frac{nb}{2r\alpha}}
\||f|^{\frac1{b+1}}\|^{\theta(b+1)}_{L^\infty(I;L^r)}
\||f|^{\frac1{b+1}}\|^{(1-\theta)(b+1)}_{\mc{C}_{q}(I;L^p)}
\enn
for $p\ge r(b+1),$ where $\theta$ is the same as in (i).
\end{lem}



\section{Well-posedness in Lebesgue spaces}
\setcounter{equation}{0}

In this section we consider the following Cauchy problem for the semi-linear fractional power
dissipative equation
\begin{eqnarray}\label{3.2}
\begin{cases}
u_t+(-\tr)^\alpha u=\pm |u|^bu,\ &(t,x)\in [0,\infty)\times\mr^n;\\
u(0)=\varphi(x),\ & x\in \mr^n.
\end{cases}
\end{eqnarray}
We shall study the well-posedness of the Cauchy problem (\ref{3.2}) for the initial data
$\varphi\in L^r(\mr^n)$, $\ds r\ge r_0=\frac{nb}{2\alpha}>1$.
The corresponding integral equation is
\be\label{integral}
u(x,t)=S_{\alpha}(t)\varphi(x)+\int^t_0S_{\alpha}(t-\tau)f(u(\tau,x))\md\tau
=S_\alpha(t)\varphi(x)+\mb{G}f(u)\triangleq\mc{T}(u),
\ee
where $f(u)=\pm|u|^bu$. The solution to the integral equation (\ref{integral}) is called a mild
solution which, by the standard regularity effect, is regular for $t>0$.

We first consider the solution to (\ref{3.2}) (or equivalently (\ref{integral})) in the space
\be\label{4.2}
X(I)=C(I;L^r(\mr^n))\cap L^q(I;L^p(\mr^n)),
\ee
where $I=[0,T)$ for $T>0.$
Using Lemmas \ref{lem3.2} and \ref{lem3.3} and applying the Banach contraction mapping
principle to the integral operator $\mc{T}$, it is easy to establish the following theorems
on the existence of local solutions or global small solutions to the problem (\ref{3.2}).
We omit the proof here for succinctness.

\begin{thm}\label{thm4.1}
Let $1<r_0={nb}/({2\alpha})\le r$ and let $\varphi\in L^r(\mr^n)$.
Assume that $(q,p,r)$ is an arbitrary admissible triplet.

(i) There exist $T>0$ and a unique mild solution $u\in X(I)$ to the problem (\ref{3.2}),
where $T=T(\|\varphi\|_{L^r})$ depends on the norm $\|\varphi\|_{L^r}$ for $r>r_0$, and
$T=T(\varphi)$ depends on $\varphi$ itself for the case $r=r_0.$

(ii) If $r=r_0$, then $T=\infty$ provided that $\|\varphi\|_{L^r}$ is sufficiently small.
In other words, there exists a global small solution
$u\in C_b([0,\infty);L^r(\mr^n))\cap L^q([0,\infty);L^p(\mr^n))$.

(iii) Let $[0,T^*)$ be the maximal existence interval of the solution $u$ to the problem (\ref{3.2})
(or equivalently (\ref{integral})) such that $u\in L^q([0,T^*);L^p(\mr^n))\cap C_b([0,T^*);L^r(\mr^n))$
for $r>r_0$. Then
\ben
\|u(s)\|_{L^r}\ge \frac C{(T^*-s)^{\frac1b-\frac n{2r\alpha}}}.
\enn
\end{thm}

We now consider the solution to (\ref{3.2}) (or equivalently (\ref{integral})) in the space
\be\label{4.4}
Y(I)=C_b(I;L^r(\mr^n))\cap \mc{\dot{C}}_q(I;L^p(\mr^n)),
\ee
where $I=[0,T)$ for $T>0.$ Making use of Lemmas \ref{lem3.2} and \ref{lem3.4} together with
the Banach contraction mapping principle to the integral equation (\ref{integral})
we can derive the following well-posedness results.

\begin{thm}\label{thm4.2}
Let $1<r_0={nb}/({2\alpha})\le r$ and let $\varphi\in L^r(\mr^n)$.
Assume that $(q,p,r)$ is any generalized admissible triplet.

(i) There exist $T>0$ and a unique mild solution $u\in Y(I)$ to the problem (\ref{3.2}),
where $T=T(\|\varphi\|_{L^r})$ depends on the norm $\|\varphi\|_{L^r}$ for the case $r>r_0$,
and $T=T(\varphi)$ depends on $\varphi$ itself for the case $r=r_0$.

(ii) If $r=r_0$, the $T=\infty$ provided that $\|\varphi\|_{L^{r_0}}$ is sufficiently small.
In other words, there exists a global small solution $u\in C_b([0,\infty);L^r(\mr^n))
\cap\mc{\dot{C}}_q([0,\infty);L^p(\mr^n))$.

(iii) Let $I=[0,T^*)$ be the maximal existence interval of the solution $u$ to the problem (\ref{integral})
such that $u\in C_b(I;L^r(\mr^n))\cap \mc{\dot{C}}_q(I;L^p(\mr^n))$ for $r>r_0$. Then
\ben
\|u(s)\|_{L^r}\ge \frac C{(T^*-s)^{\frac1b-\frac n{2r\alpha}}}.
\enn
\end{thm}

Our method is also valid for the case of convective nonlinear term, that is, the following Cauchy
problem for the fractional power dissipative convective equation:
\begin{eqnarray}\label{4.6}
\begin{cases}
u_t+(-\tr)^\alpha u=(a\cdot\nabla)g(u),\ &(t,x)\in [0,\infty)\times\mr^n;\\
u(0)=\varphi(x),\ & x\in \mr^n,
\end{cases}
\end{eqnarray}
where $b>0$, $\alpha>0$ and $a\in\mr^n$ is a given $n-$dimensional vector.
By Duhamel's principle the problem (\ref{4.6}) is equivalent to the integral equation:
\be\label{4.7}
u(t,x)=S_\alpha(t)\varphi(x)+\int^t_0S_\alpha(t-\tau)(a\cdot\nabla)g(u)\md\tau
\triangleq S_\alpha(t)\varphi(x)+\widetilde{\mb{G}}g(u),
\ee
where $g(u)=\pm |u|^bu$.

Arguing similarly as in the proof of Lemmas \ref{lem3.3} and \ref{lem3.4}, we
have the following nonlinear estimates.

\begin{lem}\label{lem4.1}
For $b>0$, $\alpha>1/2$ and $T>0$, let $\ds r_1=\frac{nb}{2\alpha-1}$ and $I=[0,T)$.
Assume that $r\ge r_1>1$. Let $(q,p,r)$ be an arbitrary admissible triplet satisfying that $p>b+1$.
If $\ds f\in L^\frac q{b+1}(I;L^{\frac p{b+1}}(\mr^n))$, then
\ben\nonumber
\|\widetilde{\mb{G}}f\|_{L^\infty(I;L^r)}+\|\widetilde{\mb{G}}f\|_{L^q(I;L^p)}
\le CT^{1-\frac1{2\alpha}-\frac{nb}{2r\alpha}}\|f\|_{L^{\frac q{b+1}}(I;L^{\frac p{b+1}})}
\enn
for $p<r(1+b)$, and
\ben\nonumber
\|\widetilde{\mb{G}}f\|_{L^\infty(I;L^r)}+\|\widetilde{\mb{G}}f\|_{L^q(I;L^p)}
\le CT^{1-\frac1{2\alpha}-\frac{nb}{2r\alpha}}\||f|^{\frac1{b+1}}\|^{\theta(b+1)}_{L^\infty(I;L^r)}
\||f|^{\frac1{b+1}}\|^{(1-\theta)(b+1)}_{L^q(I;L^p)}
\enn
for $p\ge r(b+1)$, where $\ds\theta=\frac{p-r(b+1)}{(b+1)(p-r)}$.
\end{lem}

\begin{lem}\label{lem4.2}
For $b>0$, $\alpha>1/2$ and $T>0$, let $\ds r_1=\frac{nb}{2\alpha-1}$ and $I=[0,T)$.
Assume that $r\ge r_1>1$. Let $(q,p,r)$ be an arbitrary generalized admissible triplet satisfying that
$p>b+1$. If $\ds f\in\mc{C}_{\frac q{b+1}}(I;L^{\frac p{b+1}}(\mr^n))$, then
\ben
\|\widetilde{\mb{G}}f\|_{L^\infty(I;L^r)}+\|\widetilde{\mb{G}}f\|_{\mc{C}_q(I;L^p)}
\le CT^{1-\frac1{2\alpha}-\frac{nb}{2r\alpha}}\|f\|_{\mc{C}_{\frac q{b+1}}(I;L^{\frac p{b+1}})}
\enn
for $p< r(1+b)$, and
\ben\nonumber
\|\widetilde{\mb{G}}f\|_{L^\infty(I;L^r)}+\|\widetilde{\mb{G}}f\|_{\mc{C}_q(I;L^p)}
\le CT^{1-\frac1{2\alpha}-\frac{nb}{2r\alpha}}\||f|^{\frac1{b+1}}\|^{\theta(b+1)}_{L^\infty(I;L^r)}
\||f|^{\frac1{b+1}}\|^{(1-\theta)(b+1)}_{\mc{C}_{q}(I;L^p)}
\enn
for $p\ge r(b+1)$, where $\ds\theta=\frac{p-r(b+1)}{(b+1)(p-r)}$.
\end{lem}

Using Lemmas \ref{lem3.2} and \ref{lem4.1} together with the Banach contraction mapping principle
we can get the well-posedness in the space $X(I)$ defined by (\ref{4.2}) of the Cauchy problem (\ref{4.6}).

\begin{thm}\label{thm4.3}
Let $\ds 1<r_1=\frac{nb}{2\alpha-1}\le r,$ $\varphi\in L^r(\mr^n)$ and for $T>0$ let $I=[0,T).$
Assume that $(q,p,r)$ is an arbitrary admissible triplet.

(i) There exist a $T>0$ and a unique mild solution to the problem (\ref{4.6}) such that
$u\in X(I)$, where $T=T(\|\varphi\|_{L^r})$ depends on the norm $\|\varphi\|_{L^r}$ for
the case $r>r_1$, or $T=T(\varphi)$ depends on $\varphi$ itself for the case $r=r_1.$

(ii) If $r=r_1$, then we can take $T=\infty$ provided that $\|\varphi\|_{L^r}$ is sufficiently
small. In other words, there exists a global small solution
$u\in C_b([0,\infty);L^r(\mr^n))\cap L^q((0,\infty);L^p(\mr^n))$.

(iii) Let $I=[0,T^*)$ be the maximal existence interval of the solution $u$ to the problem
(\ref{4.6}) (or equivalently (\ref{4.7})) such that
$u\in L^q(I;L^p(\mr^n))\cap C_b(I;L^r(\mr^n))$ for $r>r_1$. Then
\ben
\|u(s)\|_{L^r}\ge \frac C{(T^*-s)^{\frac1b-\frac1{2b\alpha}-\frac n{2r\alpha}}}.
\enn
\end{thm}

Similarly, making use of Lemmas \ref{lem3.2} and \ref{lem4.2} and the Banach contraction
mapping principle we can establish the following well-posedness in the space $Y(I)$ defined
in (\ref{4.4}) of the Cauchy problem (\ref{4.6}).

\begin{thm}
Let $\ds 1<r_1=\frac{nb}{2\alpha-1}\le r$, $\varphi\in L^r(\mr^n)$ and for $T>0$ let $I=[0,T).$
Assume that $(q,p,r)$ is an arbitrary generalized admissible triplet.

(i) There exist a $T>0$ and a unique mild solution to the problem (\ref{4.6}) such that
$u\in Y(I)$, where $T=T(\|\varphi\|_{L^r})$ depends on the norm $\|\varphi\|_{L^r}$ for the case
$r>r_1$ or $T=T(\varphi)$ depends on $\varphi$ itself for the case $r=r_1$.

(ii) If $r=r_1$, then $T=\infty$ provided that $\|\varphi\|_{L^r}$ is sufficiently small, that is,
there exists a global small solution
$u\in C_b([0,\infty);L^r(\mr^n))\cap\mc{\dot{C}}_q([0,\infty);L^p(\mr^n))$.

(iii) Let $I=[0,T^*)$ be the maximal existence interval of the solution $u$ to the
problem (\ref{4.6}) such that $u\in C_b(I;L^r(\mr^n))\cap \mc{\dot{C}}_q(I;L^p(\mr^n))$
for $r>r_1$. Then
\ben
\|u(s)\|_{L^r}\ge \frac C{(T^*-s)^{\frac1b-\frac1{2b\alpha}-\frac n{2r\alpha}}}.
\enn
\end{thm}

\section{Fractional power dissipative equations with more general nonlinear terms}
\subsection{The case of more general nonlinear terms}
\setcounter{equation}{0}

In this subsection we study well-posedness in Lebesgue spaces for the case of more general nonlinear terms.
In particular, we consider the following cases:
\be\label{5.1}
\begin{cases}
u_t+(-\tr)^{\alpha}u=f_1(u)+f_2(u),\ &(t,x)\in [0,\infty)\times \mr^n,\ \alpha>0,\\
u(0,x)=\varphi(x),\ & x\in \mr^n;
\end{cases}
\ee
\be\label{5.2}
\begin{cases}
u_t+(-\tr)^{\alpha}u=f_1(u)+(\beta\cdot\nabla)f_2(u),\ &(t,x)\in[0,\infty)\times\mr^n,\ 2\alpha>1,\\
u(0,x)=\varphi(x),\ & x\in\mr^n;
\end{cases}
\ee
and
\be\label{5.3}
\begin{cases}
u_t+(-\tr)^{\alpha}u=f_1(u)+\nabla^2f_2(u),\ &(t,x)\in[0,\infty)\times\mr^n,\ \alpha>1,\\
u(0,x)=\varphi(x),\ & x\in\mr^n.
\end{cases}
\ee
Here $f_1(u)=\pm |u|^{b_1}u,$ $f_2(u)=\pm |u|^{b_2}u$ and $\beta\in\mr^n$.
Without loss of generality we assume $b_1>b_2>0$.
Set $r_0={nb_1}/({2\alpha})$ and let $(q,p,r)$ be an arbitrary admissible or generalized admissible
triplet for $r\ge r_0$. For $T>0$ let $I=[0,T)$ and let
\ben
X(I)&=&C(I;L^r(\mr^n))\cap L^q(I;L^p(\mr^n)),\\
Y(I)&=&C(I;L^r(\mr^n))\cap \mc{C}_q(I;L^p(\mr^n)).
\enn
Then, similarly to Lemmas \ref{lem3.3} and \ref{lem3.4} we have the following variant space-time
estimates for the operator $\mb{G}$ (cf. (\ref{3.3}) for its definition).

\begin{lem}\label{lem5.1}
Assume that $r\ge r_0>1$. Let $(q,p,r)$ be an arbitrary admissible triplet satisfying that $p>b_1+1$.
If $\ds f_1\in L^\frac q{b_1+1}(I;L^{\frac p{b_1+1}}(\mr^n))$ and
$\ds f_2\in L^\frac q{b_2+1}(I;L^{\frac p{b_2+1}}(\mr^n))$, then
\ben\nonumber
&&\|\mb{G}(f_1+f_2)\|_{L^\infty(I;L^r)}+\|\mb{G}(f_1+f_2)\|_{L^q(I;L^p)}\\
&&\qquad\le CT^{1-\frac{nb_1}{2r\alpha}}\|f_1\|_{L^{\frac q{b_1+1}}(I;L^{\frac p{b_1+1}})}
+CT^{1-\frac{nb_2}{2r\alpha}}\|f_2\|_{L^{\frac q{b_2+1}}(I;L^{\frac p{b_2+1}})}
\enn
for the case $p<r(1+b_2),$ and
\ben\nonumber
&&\|\mb{G}(f_1+f_2)\|_{L^\infty(I;L^r)}+\|\mb{G}(f_1+f_2)\|_{L^q(I;L^p)}\\
&&\qquad\le CT^{1-\frac{nb_1}{2r\alpha}}\||f_1|^{\frac1{b_1+1}}\|^{\theta_1(b_1+1)}_{L^\infty(I;L^r)}
\||f_1|^{\frac1{b_1+1}}\|^{(1-\theta_1)(b_1+1)}_{L^q(I;L^p)}\\
&&\qquad\;+CT^{1-\frac{nb_2}{2r\alpha}}\||f_2|^{\frac1{b_2+1}}\|^{\theta_2(b_2+1)}_{L^\infty(I;L^r)}
\||f_2|^{\frac1{b_2+1}}\|^{(1-\theta_2)(b_2+1)}_{L^q(I;L^p)}
\enn
for the case $p\ge r(b_1+1)$, where $\ds\theta_1=\frac{p-r(b_1+1)}{(b_1+1)(p-r)}$ and
$\ds\theta_2=\frac{p-r(b_2+1)}{(b_2+1)(p-r)}.$
If $r(1+b_2)\le p<r(b_1+1)$, then
\ben\nonumber
&&\|\mb{G}(f_1+f_2)\|_{L^\infty(I;L^r)}+\|\mb{G}(f_1+f_2)\|_{L^q(I;L^p)}\\ \nonumber
&&\quad\le CT^{1-\frac{nb_1}{2r\alpha}}\|f_1\|_{L^{\frac q{b_1+1}}(I;L^{\frac p{b_1+1}})}
+CT^{1-\frac{nb_2}{2r\alpha}}\||f_2|^{\frac1{b_2+1}}\|^{\theta_2(b_2+1)}_{L^\infty(I;L^r)}
\||f_2|^{\frac1{b_2+1}}\|^{(1-\theta_2)(b_2+1)}_{L^q(I;L^p)}.
\enn
\end{lem}


\begin{lem}\label{lem5.2}
Assume that $r\ge r_0>1$. Let $(q,p,r)$ be an arbitrary generalized admissible triplet
satisfying that $p>b_1+1$.
If $\ds f_1\in \mc{C}_{\frac q{b_1+1}}(I;L^{\frac p{b_1+1}}(\mr^n))$ and
$\ds f_2\in \mc{C}_{\frac q{b_2+1}}(I;L^{\frac p{b_2+1}}(\mr^n))$, then
\ben\nonumber
&&\|\mb{G}(f_1+f_2)\|_{L^\infty(I;L^r)}+\|\mb{G}(f_1+f_2)\|_{\mc{C}_q(I;L^p)}\\
&&\qquad\le CT^{1-\frac{nb_1}{2r\alpha}}\|f_1\|_{\mc{C}_{\frac q{b_1+1}}(I;L^{\frac p{b_1+1}})}
+CT^{1-\frac{nb_2}{2r\alpha}}\|f_2\|_{\mc{C}_{\frac q{b_2+1}}(I;L^{\frac p{b_2+1}})}
\enn
for the case $p< r(1+b_2)$, and
\ben\nonumber
&&\|\mb{G}(f_1+f_2)\|_{L^\infty(I;L^r)}+\|\mb{G}(f_1+f_2)\|_{\mc{C}_q(I;L^p)}\\ \nonumber
&&\qquad\le CT^{1-\frac{nb_1}{2r\alpha}}\||f_1|^{\frac1{b_1+1}}\|^{\theta_1(b_1+1)}_{L^\infty(I;L^r)}
\||f_1|^{\frac1{b_1+1}}\|^{(1-\theta_1)(b_1+1)}_{\mc{C}_{q}(I;L^p)}\\
&&\qquad\;+CT^{1-\frac{nb_2}{2r\alpha}}\||f_2|^{\frac1{b_2+1}}\|^{\theta_2(b_2+1)}_{L^\infty(I;L^r)}
\||f_2|^{\frac1{b_2+1}}\|^{(1-\theta_2)(b_2+1)}_{\mc{C}_{q}(I;L^p)}
\enn
for the case $p\ge r(b_1+1)$, where $\theta_1$ and $\theta_2$ are the same as defined in Lemma \ref{lem5.1}.
If $r(b_2+1)\le p<r(b_1+1)$, then
\ben\nonumber
&&\|\mb{G}(f_1+f_2)\|_{L^\infty(I;L^r)}+\|\mb{G}(f_1+f_2)\|_{\mc{C}_q(I;L^p)}\\ \nonumber
&&\quad\le CT^{1-\frac{nb_1}{2r\alpha}}\|f_1\|_{\mc{C}_{\frac q{b_1+1}}(I;L^{\frac p{b_1+1}})}
+CT^{1-\frac{nb_2}{2r\alpha}}\||f_2|^{\frac1{b_2+1}}\|^{\theta_2(b_2+1)}_{L^\infty(I;L^r)}
\||f_2|^{\frac1{b_2+1}}\|^{(1-\theta_2)(b_2+1)}_{\mc{C}_{q}(I;L^p)}.
\enn
\end{lem}

Using Lemmas \ref{lem3.2} and \ref{lem5.1} and the space $X(I)$ it is easy to prove the
well-posedness of the Cauchy problem (\ref{5.1}) by the Banach contraction mapping principle.

\begin{thm}\label{thm5.1}
For $r\ge r_0>1$ let $\varphi\in L^r(\mr^n)$. Let $(q,p,r)$ be an arbitrary admissible triplet.

(i) There exist a $T>0$ and a unique mild solution to the problem (\ref{5.1})
such that $u\in X(I)$, where $T=T(\|\varphi\|_{L^r})$ depends on the norm $\|\varphi\|_{L^r}$
for the case $r>r_0$ or $T=T(\varphi)$ depends on $\varphi$ itself for the case $r=r_0.$

(ii) Let $I=[0,T^*)$ be the maximal existence interval of the solution $u$ to the problem (\ref{5.1})
such that $u\in L^q([0,T^*);L^p(\mr^n))\cap C_b([0,T^*);L^r(\mr^n))$ for $r>r_0$. Then
\ben
\|u(s)\|_{L^r}\ge \frac C{(T^*-s)^{\frac1{b_1}-\frac n{2r\alpha}}}.
\enn
Consequently, if $T^*<\infty$, then
\ben
\lim_{t\rightarrow T^*}\|u(t)\|_{L^r}=\infty.
\enn
\end{thm}

Similarly, utilizing the space $Y(I)$ and Lemmas \ref{lem3.2} and \ref{lem5.2} in conjunction
with the Banach contraction mapping principle we can establish the following
well-posedness of the Cauchy problem \ref{5.1}.

\begin{thm}\label{thm5.2}
For $r\ge r_0>1$ let $\varphi\in L^r(\mr^n)$. Let $(q,p,r)$ be an arbitrary generalized
admissible triplet.

(i) There exist a $T>0$ and a unique mild solution to the problem (\ref{5.1})
such that $u\in Y(I)$, where $T=T(\|\varphi\|_{L^r})$ depends on the norm $\|\varphi\|_{L^r}$
for the case $r>r_0$ or $T=T(\varphi)$ depends on $\varphi$ itself for the case $r=r_0.$

(ii) Let $I=[0,T^*)$ be the maximal existence interval of the solution $u$ to the problem (\ref{5.1})
such that $u\in\mc{\dot{C}}_q([0,T^*);L^p(\mr^n))\cap C_b([0,T^*);L^r(\mr^n))$ for $r>r_0$. Then
\ben
\|u(s)\|_{L^r}\ge \frac C{(T^*-s)^{\frac1{b_1}-\frac n{2r\alpha}}}.
\enn
Consequently, if $T^*<\infty$, then
\ben
\lim_{t\rightarrow T^*}\|u(t)\|_{L^r}=\infty.
\enn
\end{thm}

Consider now the problems (\ref{5.2}) and (\ref{5.3}) with the convective effect or with higher-order
derivative term, respectively, and for $0\le d<2\alpha$ let $r_d={nb_1}/({2\alpha-d}).$
Similarly to Lemmas \ref{lem5.1} and \ref{lem5.2}, we have the following estimates of the nonlinear terms.

\begin{lem}\label{lem5.3}
Let $r\ge r_d>1$ and let $(q,p,r)$ be an arbitrary admissible triplet satisfying that $p>b_1+1$.
If $\ds f_1\in L^\frac q{b_1+1}(I;L^{\frac p{b_1+1}}(\mr^n))$
and $\ds f_2\in L^\frac q{b_2+1}(I;L^{\frac p{b_2+1}}(\mr^n)),$ then
\ben\nonumber
&&\|\mb{G}(f_1+g)\|_{L^\infty(I;L^r)}+\|\mb{G}(f_1+g)\|_{L^q(I;L^p)}\\
&&\quad\le CT^{1-\frac{nb_1}{2r\alpha}}\|f_1\|_{L^{\frac q{b_1+1}}(I;L^{\frac p{b_1+1}})}
+CT^{1-\frac d{2\alpha}-\frac{nb_2}{2r\alpha}}\|f_2\|_{L^{\frac q{b_2+1}}(I;L^{\frac p{b_2+1}})}
\enn
for $p<r(1+b_2),$ and
\ben\nonumber
&&\|\mb{G}(f_1+g)\|_{L^\infty(I;L^r)}+\|\mb{G}(f_1+g)\|_{L^q(I;L^p)}\\ \nonumber
&&\quad\le CT^{1-\frac{nb_1}{2r\alpha}}\||f_1|^{\frac1{b_1+1}}\|^{\theta_1(b_1+1)}_{L^\infty(I;L^r)}
\||f_1|^{\frac1{b_1+1}}\|^{(1-\theta_1)(b_1+1)}_{L^q(I;L^p)}\\
&&\qquad+ CT^{1-\frac d{2\alpha}-\frac{nb_2}{2r\alpha}}
\||f_2|^{\frac1{b_2+1}}\|^{\theta_2(b_2+1)}_{L^\infty(I;L^r)}
\||f_2|^{\frac1{b_2+1}}\|^{(1-\theta_2)(b_2+1)}_{L^q(I;L^p)}
\enn
for $p\ge r(b_1+1)$, where $\theta_1$ and $\theta_2$ are the same as in Lemma \ref{lem5.1}.
If $r(1+b_2)\le p<r(b_1+1)$, then
\ben\nonumber
&&\|\mb{G}(f_1+g)\|_{L^\infty(I;L^r)}+\|\mb{G}(f_1+g)\|_{L^q(I;L^p)}\\
&&\quad\le CT^{1-\frac{nb_1}{2r\alpha}}\|f_1\|_{L^{\frac q{b_1+1}}(I;L^{\frac p{b_1+1}})}
+CT^{1-\frac d{2\alpha}-\frac{nb_2}{2r\alpha}}
\||f_2|^{\frac1{b_2+1}}\|^{\theta_2(b_2+1)}_{L^\infty(I;L^r)}
\||f_2|^{\frac1{b_2+1}}\|^{(1-\theta_2)(b_2+1)}_{L^q(I;L^p)}.
\enn
Here $g=(\beta\cdot\nabla)f_2,$ $d=1$ for the problem (\ref{5.2}) or
$g=\nabla^2f_2,$ $d=2$ for the problem (\ref{5.3}).
\end{lem}


\begin{lem}\label{lem5.4}
Let $r\ge r_d>1$ and let $(q,p,r)$ be an arbitrary generalized admissible triplet satisfying that
$p>b_1+1$. If $\ds f_1\in \mc{C}_{\frac q{b_1+1}}(I;L^{\frac p{b_1+1}}(\mr^n))$
and $\ds f_2\in \mc{C}_{\frac q{b_2+1}}(I;L^{\frac p{b_2+1}}(\mr^n))$, then
\ben\nonumber
&&\|\mb{G}(f_1+g)\|_{L^\infty(I;L^r)}+\|\mb{G}(f_1+g)\|_{\mc{C}_q(I;L^p)}\\
&&\qquad\le CT^{1-\frac{nb_1}{2r\alpha}}\|f_1\|_{\mc{C}_{\frac q{b_1+1}}(I;L^{\frac p{b_1+1}})}
+CT^{1-\frac d{2\alpha}-\frac{nb_2}{2r\alpha}}\|f_2\|_{\mc{C}_{\frac q{b_2+1}}(I;L^{\frac p{b_2+1}})}
\enn
for $p<r(1+b_2),$ and
\ben\nonumber
&&\|\mb{G}(f_1+g)\|_{L^\infty(I;L^r)}+\|\mb{G}(f_1+g)\|_{\mc{C}_q(I;L^p)}\\
&&\quad\le CT^{1-\frac{nb_1}{2r\alpha}}\||f_1|^{\frac1{b_1+1}}\|^{\theta_1(b_1+1)}_{L^\infty(I;L^r)}
\||f_1|^{\frac1{b_1+1}}\|^{(1-\theta_1)(b_1+1)}_{\mc{C}_{q}(I;L^p)}\\
&&\qquad +CT^{1-\frac d{2\alpha}-\frac{nb_2}{2r\alpha}}
\||f_2|^{\frac1{b_2+1}}\|^{\theta_2(b_2+1)}_{L^\infty(I;L^r)}
\||f_2|^{\frac1{b_2+1}}\|^{(1-\theta_2)(b_2+1)}_{\mc{C}_{q}(I;L^p)}
\enn
for $p\ge r(b_1+1)$, where $\theta_1$ and $\theta_2$ are the same as in Lemma \ref{lem5.1}.
If $r(b_2+1)\le p<r(b_1+1)$, then
\ben\nonumber
&&\|\mb{G}(f_1+g)\|_{L^\infty(I;L^r)}+\|\mb{G}(f_1+g)\|_{\mc{C}_q(I;L^p)}\\
&&\quad\le CT^{1-\frac{nb_1}{2r\alpha}}\|f_1\|_{\mc{C}_{\frac q{b_1+1}}(I;L^{\frac p{b_1+1}})}
+ CT^{1-\frac d{2\alpha}-\frac{nb_2}{2r\alpha}}
\||f_2|^{\frac1{b_2+1}}\|^{\theta_2(b_2+1)}_{L^\infty(I;L^r)}
\||f_2|^{\frac1{b_2+1}}\|^{(1-\theta_2)(b_2+1)}_{\mc{C}_{q}(I;L^p)}.
\enn
Here $g$ and $d$ are the same as defined in Lemma \ref{lem5.3}.
\end{lem}

Similarly as before, using Lemma \ref{lem3.2} and Lemma \ref{lem5.3} or \ref{lem5.4}
and the space $X(I)$ or $Y(I)$ we can establish the following results
(Theorem \ref{thm5.3} or \ref{thm5.4}, respectively)
on well-posedness of the Cauchy problem (\ref{5.2}) and (\ref{5.3}).

\begin{thm}\label{thm5.3}
For $r\ge r_d={nb_1}/({2\alpha-d})>1$ with $0\le d<2\alpha$ let $\varphi\in L^r(\mr^n)$.
Assume that $(q,p,r)$ is an arbitrary admissible triplet.

(i) There exist a $T>0$ and a unique mild solution $u\in X(I)$ to the problem
(\ref{5.2}) or (\ref{5.3}), where $T=T(\|\varphi\|_{L^r})$ depends on the
norm $\|\varphi\|_{L^r}$ for the case $r>r_d$ or $T=T(\varphi)$ depends
on $\varphi$ itself for the case $r=r_d.$

(ii) Let $I=[0,T^*)$ be the maximal existence interval of the solution $u$ to the problem
(\ref{5.2}) or (\ref{5.3}) such that $u\in L^q([0,T^*);L^p(\mr^n))\cap C_b([0,T^*);L^r(\mr^n))$
for $r>r_d$. Then
\ben
\|u(s)\|_{L^r}\ge \frac C{(T^*-s)^{\frac1{b_1}-\frac d{2b_1\alpha}-\frac n{2r\alpha}}}.
\enn
Consequently, if $T^*<\infty$, then
\ben
\lim_{t\rightarrow T^*}\|u(t)\|_{L^r}=\infty.
\enn
Here $d=1$ in the case of (\ref{5.2}) or $d=2$ in the case of (\ref{5.3}).
\end{thm}

\begin{thm}\label{thm5.4}
For $r\ge r_d={nb_1}/({2\alpha-d})>1$ with $0\le d<2\alpha$ let $\varphi\in L^r(\mr^n)$.
Assume that $(q,p,r)$ is an arbitrary generalized admissible triplet.

(i) There exist a $T>0$ and a unique mild solution $u\in X(I)$ to the problem
(\ref{5.2}) or (\ref{5.3}), where $T=T(\|\varphi\|_{L^r})$ depends on the
norm $\|\varphi\|_{L^r}$ for the case $r>r_d$ or $T=T(\varphi)$ depends
on $\varphi$ itself for the case $r=r_d.$

(ii) Let $[0,T^*)$ be the maximal existence interval of the solution $u$ to (\ref{5.2})
or (\ref{5.3}) such that $u\in\mc{C}_q([0,T^*);L^p(\mr^n))\cap C_b([0,T^*);L^r(\mr^n))$ for
$r>r_d$. Then
\ben
\|u(s)\|_{L^r}\ge \frac C{(T^*-s)^{\frac1{b_1}-\frac d{2b_1\alpha}-\frac n{2r\alpha}}}.
\enn
Consequently, if $T^*<\infty$, then
\ben
\lim_{t\rightarrow T^*}\|u(t)\|_{L^r}=\infty.
\enn
Here $d=1$ in the case of (\ref{5.2}) or $d=2$ in the case of (\ref{5.3}).
\end{thm}

\begin{rem}\label{r5.1} {\rm
(1) Our methods can also be applied to the case of general nonlinear terms
\ben
f(u)=\sum_{k=1}^mP_k(D)f_k(u),
\enn
where $P_k(D)$ is a homogeneous pseudo-differential operator of order $d_k\in [0,2\alpha)$
and $f_k(u)$ behaves like $|u|^{b_k}u$ or $|u|^{b_k+1}$. In fact, we only need to take
$\ds r_d=\max\limits_{1\le k\le m}\frac{nb_k}{2\alpha-d_k},$ and the similar results to Theorems
\ref{thm5.3} and \ref{thm5.4} still hold.

(2) One easily sees by the scaling argument that the solution space corresponds to the
`highest' nonlinear growth, that is, $r\ge r_d>1$.
}
\end{rem}

Theorems \ref{thm5.3} and \ref{thm5.4} can be applied to get the local well-posedness
of the following Ginzburg-Landau equations for the population model \cite{CM}:
\ben
\begin{cases}
u_t+a_1\nabla^4 u=G(u)+a_2\nabla^2u+a\nabla^2 u^3,\ &(t,x)\in [0,\infty)\times \mr^n,\\
u(0,x)=\varphi(x),\ & x\in \mr^n,
\end{cases}
\enn
or the generalized Ginzburg-Landau equations \cite{Chen}:
\ben
\begin{cases}
u_t+a_1\nabla^4 u=G(u)+a_2\nabla^2u+\nabla^2 g(u),\ &(t,x)\in [0,\infty)\times \mr^n,\\
u(0,x)=\varphi(x),\ & x\in \mr^n,
\end{cases}
\enn
where $a_1>0$, $a>0$ and $a_2\not=0$, $G(u)$ and $g(u)$ can be of polynomial growth such as
$G(u)=c_ku^k+c_{k-1}u^{k-1}+\cdots+c_1 u$ and $g(u)=d_lu^l+d_{l-1}u^{l-1}+\cdots+d_1 u$ for $k,\,l>1$.

\subsection{Global existence for small initial data}

In this subsection we consider the problem (\ref{1.1}) with the nonlinear term
$F(u)=f_1(u)+f_2(u)$ or $F(u)=Q(D)(f_1(u)+f_2(u))$, where $f_1(u)=\pm|u|^{b_1}u,$
$f_2(u)=\pm|u|^{b_2}u$ or $f_1(u)=\pm|u|^{b_1+1},$ $f_2(u)=\pm|u|^{b_2+1}$ for $b_1\ge b_2>0$
and $Q(D)$ is a homogeneous pseudo-differential operator of order $d\in [0,2\alpha)$.
Let $r_j={nb_j}/({2\alpha-d})$ for $j=1,2.$ We show that, if the initial data
$\varphi\in L^{r_1}\cap L^{r_2}$ and $\|\varphi;L^{r_1}\cap L^{r_2}\|$ is small enough
then the problem (\ref{1.1}) has a unique global solution, as stated in the following theorem.

\begin{thm}\label{thm5.5}
Let $r_j>1$ and $(q_j, p_j, r_j)$ be generalized admissible triplets such that
$p_j\le r_j(1+b_j)$ for $j=1,2$. Then there exists a unique solution
$u\in C([0,\infty);L^{r_1}\cap L^{r_2})\cap\mc{\dot{C}}_{q_1}([0,\infty);L^{p_1})\cap
\mc{\dot{C}}_{q_2}([0,\infty);L^{p_2})$ to the Cauchy problem (\ref{1.1})
provided that the initial data $\varphi\in L^{r_1}\cap L^{r_2}$ and
the norm $\|\varphi;L^{r_1}\cap L^{r_2}\|<\delta$ for some sufficiently small $\delta>0$.
In particular, if $b_1=b_2$ we recover the global existence of small solutions to the Cauchy problem
(\ref{3.2}).
\end{thm}

\begin{proof}
The proof is broken down into the following three steps.

Step 1. Assume that the generalized admissible triplets $(q_j,p_j,r_j)$
($j=i,1$) satisfy the conditions
\be\label{condition1}
1+b_j<p_j\le r_j(1+b_j),\;\; j=1,\,2
\ee
and
\be\label{condition2}
\frac {r_1}{p_1}=\frac{r_2}{p_2}.
\ee
For $I=[0,\infty)$ define the solution space as
\ben
Z(I)=\{u\,\big|\,u\in C_b(I;L^{r_1}\cap L^{r_2})\cap\mc{\dot{C}}_{q_1}(I;L^{p_1})
\cap\mc{\dot{C}}_{q_2}(I;L^{p_2})\}
\enn
with the norm
\ben
\|u;Z(I)\|=\sum_{j=1}^2\sup_{t\in I}t^{\frac1{q_j}}\|u\|_{L^{p_j}}
+\sum_{j=1}^2\sup_{t\in I}\|u\|_{L^{r_j}}.
\enn
The problem (\ref{1.1}) can be written in the integral form as
\ben
u(x,t)=S_{\alpha}(t)\varphi(x)+\mb{\tilde{G}}(f_1(u)+f_2(u))\triangleq\mc{T}u,
\enn
where
$$
\mb{\tilde{G}}(f_1(u)+f_2(u))=\int^t_0S_{\alpha}(t-\tau)Q(D)[f_1(u(\tau,x))+f_2(u(\tau,x))]\md\tau.
$$

Now consider the operator $\mc{T}$ in the complete metric space
\ben
E(I)=\{u\in Z(I)\,\big|\,\|u;Z(I)\|\le \delta\}
\enn
with the metric
\ben
d(u,v)=\|u-v; Z(I)\|,\;\; u,v\in E(I),
\enn
where $\delta>0$ is a sufficiently small constant to be determined later.
By Lemmas \ref{lem3.1} and \ref{lem3.2} one has
\be\label{5.36}
\|S_\alpha(t)\varphi;Z(I)\|\le C(\|\varphi\|_{r_1}+\|\varphi\|_{r_2}).
\ee
By Lemma \ref{lem5.4} and the H\"older inequality for $\ds\frac1p=\frac1{p_1}+\frac{b_2}{p_2}$
we get
\be\nonumber
&&\|\mb{\tilde{G}}(f_1(u)+f_2(u));\mc{C}_{q_1}(I;L^{p_1})\|\\ \nonumber
&&\le C\|u;\mc{C}_{q_1}(I;L^{p_1})\|^{b_1+1}+\sup_{t\in I}t^{\frac1{q_1}}
\int^t_0(t-\tau)^{-\frac d{2\alpha}-\frac{nb_2}{2\alpha p_2}}\|u\|^{b_2}_{p_2}\|u\|_{p_1}\md\tau\\ \no
&&\le C\|u;\mc{C}_{q_1}(I;L^{p_1})\|^{b_1+1}
+C\int^1_0(1-\tau)^{-\frac d{2\alpha}-\frac{nb_2}{2\alpha p_2}}
\tau^{-\frac1{q_1}-\frac{\alpha_2}{q_2}}\md\tau\|u;\mc{C}_{q_2}(I;L^{p_2})\|^{b_2}
\|u;\mc{C}_{q_1}(I;L^{p_1})\|\\ \label{5.37}
&&\le C\sum_{j=1}^2\|u;\mc{C}_{q_j}(I;L^{p_j})\|^{b_j}\|u;\mc{C}_{q_1}(I;L^{p_1})\|
\ee
and
\be\nonumber
&&\|\mb{\tilde{G}}(f_1(u)+f_2(u));C(I;L^{r_1})\|\\ \nonumber
&&\le C\sup_{t\in I}\int^t_0(t-\tau)^{-\frac d{2\al}-\frac{n}{2\al}(\frac{1+b_1}{p_1}-\frac1{r_1})}
\|u\|_{p_1}^{1+b_1}\md\tau+C\sup_{t\in I}
\int^t_0(t-\tau)^{-\frac d{2\al}-\frac n{2\al}(\frac1p-\frac1{r_1})}\|u\|_{p_2}^{b_2}\|u\|_{p_1}\md\tau\\ \no
&&\le C\sup_{t\in I}\int^t_0(t-\tau)^{-\frac d{2\alpha}-\frac {n}{2\alpha}
(\frac{1+b_1}{p_1}-\frac1{r_1})}\md\tau\|u;\mc{C}_{q_1}(I;L^{p_1})\|^{1+b_1}\\ \nonumber
&&+C\sup_{t\in I}\int^t_0(t-\tau)^{-\frac d{2\alpha}-\frac n{2\alpha}(\frac1p-\frac1{r_1})}
\tau^{-\frac{b_2}{q_2}-\frac1{q_1}}\md\tau
\|u;\mc{C}_{q_2}(I;L^{p_2})\|^{b_2}\|u;\mc{C}_{q_1}(I;L^{p_1})\|\\ \label{5.38}
&&\le C\sum_{j=1}^2\|u;\mc{C}_{q_j}(I;L^{p_j})\|^{b_j}\|u;\mc{C}_{q_1}(I;L^{p_1})\|.
\ee
Similarly, we have
\be\label{5.39}
\|\mb{\tilde{G}}(f_1(u)+f_2(u));\mc{C}_{q_2}(I;L^{p_2})\|
&\le& C\sum_{j=1}^2\|u;\mc{C}_{q_j}(I;L^{p_j})\|^{b_j}\|u;\mc{C}_{q_2}(I;L^{p_2})\|,\\ \label{5.40}
\|\mb{\tilde{G}}(f_1(u)+f_2(u));C(I;L^{r_1})\|&\le&
C\sum_{j=1}^2\|u;\mc{C}_{q_j}(I;L^{p_j})\|^{b_j}\|u;\mc{C}_{q_2}(I;L^{p_2})\|.
\ee
Combining the estimates (\ref{5.36})-(\ref{5.40}) and choosing $\delta>0$ small enough we get
\be\nonumber
\|\mc{T}u;Z(I)\|&\le& C\|\varphi;L^{r_1}\cap L^{r_2}\|+C\sum_{j=1}^2\|u;\mc{C}_{q_j}(I;L^{p_j})
\|^{b_j}\|u;\mc{C}_{q_1}(I;L^{p_1})\|+ \\ \nonumber
&& C\sum_{j=1}^2\|u;\mc{C}_{q_j}(I;L^{p_j})\|^{b_j}\|u;\mc{C}_{q_2}(I;L^{p_2})\|\\ \label{5.41}
&\le& C\|\varphi;L^{r_1}\cap L^{r_2}\|+C\delta^{b_1+1}+C\delta^{b_2+1}<\delta
\ee
provided that $C\|\varphi;L^{r_1}\cap L^{r_2}\|<\frac\delta 2$.
Noting the definition of $E(I)$ and the fact that
\ben
d(\mc{T}u,\mc{T}v)\le C\sum_{j=1}^2(\|u;\mc{C}_{q_j}(I;L^{p_j})
\|^{b_j}+\|v;\mc{C}_{q_j}(I;L^{p_j}) \|^{b_j})d(u,v),
\enn
one has $d(\mc{T}u,\mc{T}v)\le \frac12 d(u,v)$. Furthermore, from (\ref{5.41}), and since
\ben
\lim_{t\rightarrow 0^+}t^{\frac1{q_j}}\|S_{\alpha}(t)\varphi;L^{r_1}\cap L^{r_2}\|=0,
\enn
it follows that
\ben
\lim_{t\rightarrow 0^+}t^{\frac1{q_j}}\|\mc{T}u\|_{L^{p_j}}=0,\;\;j=1,2.
\enn
Thus $\mc{T}$ is a contraction mapping from $E(I)$ into itself so,
by the Banach contraction mapping principle there exists a unique solution $u\in E(I)$.

Step 2. We show that $u\in\mc{\dot{C}}_{q}(I;L^{p})$ for any generalized admissible
triplet $(q,p,r_1)$ satisfying the condition (\ref{condition1}).
Without lost of generality we assume that the generalized admissible triplet $(q_j,p_j,r_j)$
satisfies the conditions (\ref{condition1}) and (\ref{condition2}) for $j=1,2$.
Arguing similarly as in deriving (\ref{5.37}) we have
\ben\nonumber
\|u;\mc{C}_q(I;L^p)\|&\le& C\|\varphi\|_{L^{r_1}}+\sup_{t\in I}t^{\frac1q}
\int^t_0(t-\tau)^{-\frac d{2\alpha}-\frac{nb_1}{2\alpha p_1}}\|u\|^{b_1}_{p_1}\|u\|_{p}\md\tau\\
&&+\sup_{t\in I}t^{\frac1q}\int^t_0(t-\tau)^{-\frac d{2\alpha}-\frac{nb_2}{2\al p_2}}
\|u\|^{b_2}_{p_2}\|u\|_{p}\md\tau\\ \nonumber
&\le& C\|\varphi\|_{L^{r_1}}+C\|u;\mc{C}_{q_1}(I;L^{p_1})\|^{b_1}\|u;\mc{C}_{q}(I;L^{p})\|\\
&&+C\|u;\mc{C}_{q_2}(I;L^{p_2})\|^{b_2}\|u;\mc{C}_{q}(I;L^{p})\|\\
&\le& C\|\varphi\|_{L^{r_1}}+C(\delta^{b_1}+\delta^{b_2})\|u;\mc{C}_{q}(I;L^{p})\|,
\enn
which implies that for small $\delta>0$
\ben
\|u;\mc{C}_{q}(I;L^{p})\|<\infty.
\enn
Similar arguments as above give
\ben
\lim_{t\rightarrow 0^+}t^{\frac1q}\|u\|_{L^p}=0.
\enn
Thus, it is derived that $u\in\mc{\dot{C}}_{q}(I;L^{p})$ for any
generalized admissible triplet $(q,p,r_1)$ satisfying the condition (\ref{condition1}).

Step 3. Finally, if $(q,p,r_1)$ is a generalized admissible triplet satisfying $p\le 1+b_1,$
then the result that $u\in\mc{\dot{C}}_{q}(I;L^{p})$ follows by interpolation between $C(I;L^{r_1})$ and
$\mc{\dot{C}}_{\hat{q}(\hat{p},r_1)}(I;L^{\hat{p}})$ with the generalized admissible triplet
$(\hat{q},\hat{p},r_2)$ satisfying the condition (\ref{condition1}).
The proof is thus complete.
\end{proof}

\begin{rem}\label{r5.2} {\rm
One can easily see that there exist generalized admissible triplets $(q_j,p_j,r_j)$ ($j=1,2$)
satisfying condition (\ref{condition2}).
In fact, without loss of generality we assume $b_1\ge b_2$, which implies $r_1\ge r_2$.
Since
\ben
\frac{1+b_1}{r_1}<\frac{p_1}{r_1}\le (1+b_1),\qquad \frac{1+b_2}{r_2}<\frac{p_2}{r_2}\le (1+b_2)
\enn
and
\ben
\left(\frac{2\al-d}{n}\left(1+\frac1b_2\right),1+b_2\right)\subset
\left(\frac{2\al-d}{n}\left(1+\frac1b_1\right),1+b_1\right),\quad\mx{for}\; d\in[0,2\alpha),
\enn
it is not difficult to choose generalized admissible triplets $(q_j,p_j,r_j)$ ($j=1,2$)
satisfying the condition (\ref{condition2}). Moveover, since $\ds r_j=\frac{nb_j}{2\alpha-d}$,
we also have
\ben
\frac {b_1}{p_1}=\frac {b_2}{p_2},\quad \mx{and}\quad\frac{b_1}{q_1}=\frac {b_2}{q_2}.
\enn
}
\end{rem}

\section{Global well-posedness for high frequency data}
\setcounter{equation}{0}

In the previous sections we proved the well-posedness of the problem (\ref{3.2})
for the initial data $\varphi\in L^r(\mr^n)$ with $r\ge r_0={nb}/({2\al}),$ that is,
if the norm $\|\varphi\|_{L^{r_0}}$ of the initial data is small enough, then the solution $u$
exists globally. In this section we shall show that the solution $u$ exists globally
if the norm $\ds\|\varphi;\dot{B}^{\frac np-\frac{2\alpha}b}_{p,\infty}(\mr^n)\|$ is small
enough (in this case, the norm $\|\varphi\|_{L^{r_0}}$ may be large).
To this end, let $I=[0,\infty)$ and let us introduce the following solution space
\ben
X(I)=C(I;L^{r_0}(\mr^n))\cap\mc{C}_q(I;L^p(\mr^n))
\enn
with the norm
\ben
\|u\|_{X(I)}=\sup_{t>0}\|u(t)\|_{\dot{B}^{-\sigma}_{p,\infty}}+\sup_{t>0}t^{\frac1q}\|u\|_p,
\enn
where $r_0={nb}/({2\alpha})\le p<r_0(b+1)$, $p>b+1,$ $\sigma={2\al}/b-n/p\ge 0$
and $\ds\frac1q=\frac n{2\al}(\frac1{r_0}-\frac1p)$. Consider the operator $\mc{T}$,
defined in (\ref{integral}), in the complete metric space
\ben
X_\delta=\{u(t)\in X(I)\,\big|\,\|u(t)\|_{X(I)}\le 2\delta\}
\enn
with the metric $d(u,v)=\|u-v\|_{X(I)}$ for $u,v\in X_\delta$, where $\delta$ is a small constant
to be determined later.
Using the equivalent characterization of Besov spaces (see Proposition \ref{prop2.1}) we have
\be\label{6.5}
\|S_\alpha(t)\varphi\|_{X(I)}&=&\sup_{t>0}\|S_\alpha(t)\varphi\|_{\dot{B}^{-\sigma}_{p,\infty}}
     +\sup_{t>0}t^{\frac1q}\|S_\alpha(t)\varphi\|_{L^p}
\le C\|\varphi\|_{\dot{B}^{-\sigma}_{p,\infty}}.
\ee
By Lemma \ref{lem3.1} and the Sobolev embedding
$L^{r_0}(\mr^n)\hookrightarrow\dot{B}^{-\sigma}_{p,\infty}(\mr^n)$,
it is seen that
\be\nonumber
\|\mb{G}f\|_{\dot{B}^{-\sigma} _{p,\infty}}&\le&\|\mb{G}f\|_{L^{r_0}}\\ \nonumber
&\le& C\sup_{0<t<\infty}\int^t_0(t-\tau)^{-\frac n{2\al}(\frac{b+1}p-\frac{2\al}{nb})}
      \||f(u(\tau))|^{\frac1{b+1}}\|_{L^p}^{b+1}\md\tau \\ \nonumber
&\le& C\sup_{0<t<\infty}\int^t_0(t-\tau)^{-\frac n{2\al}(\frac{b+1}p-\frac{2\al}{nb})}
      \tau^{-\frac{b+1}q}\md\tau\|u\|^{b+1}_{\mc{C}_q(I;L^p)}\\ \label{6.6}
&\le& C\|u\|^{b+1}_{X(I)},\\
\|\mb{G}f\|_{\mc{C}_q(I;L^p)}&\le&\sup_{0<t<\infty}t^{\frac1q}
   \int^t_0(t-\tau)^{-\frac n{2\al}(\frac{b+1}p-\frac1p)}
   \||f(u(\tau))|^{\frac1{b+1}}\|^{b+1}_{L^p}\md\tau\\ \nonumber
&\le&\sup_{0<t<\infty}t^{\frac1q}\int^t_0(t-\tau)^{-\frac{nb}{2\alpha p}}
   \tau^{-\frac{b+1}q}\md\tau\|u\|^{b+1}_{\mc{C}_q(I;L^p)}\\ \label{6.7}
&\le& C\|u\|^{b+1}_{X(I)}.
\ee
Combining (\ref{6.5})-(\ref{6.7}) we have on noting the definition (\ref{integral}) of $\mc{T}$ that
for $u\in X_\delta$
\ben
\|\mc{T}(u)\|_{X(I)}\le C\|\varphi\|_{\dot{B}^{-\sigma} _{p,\infty}}+ C\|u\|^{b+1}_{X(I)}
\le C\|\varphi\|_{\dot{B}^{-\sigma} _{p,\infty}}+ C\delta^{b+1}.
\enn
Thus, if we take $\ds\delta=C\|\varphi\|_{\dot{B}^{-\sigma} _{p,\infty}}$ to be small enough,
then $\mc{T}$ is a contraction mapping from $X_\delta$ into itself.
The Banach contraction mapping principle implies that $\mc{T}$ has a unique fixed point
in $u\in X_\delta$ or equivalently the problem (\ref{3.2}) has a unique solution $u\in X_\delta$.
Furthermore, one can verify that
\ben\nonumber
t^{\frac1q}\|\mb{G}(u)\|_{L^p}&\le& t^{\frac1q}\int^t_0(t-\tau)^{-\frac {nb}{2\al p}}
  \tau^{-\frac{b+1}q}\md\tau(\sup_{0<\tau<t}\tau^{\frac1q}\|u(\tau)\|_{L^p})^{b+1}\\
&\le& C\delta^b\sup_{0<\tau<t}\tau^{\frac1q}\|u(\tau)\|_{L^p},
\enn
and
\ben
\lim_{t\rightarrow 0}t^{\frac1q}\|S_\alpha(t)\varphi\|_{L^p}=0.
\enn
Consequently, it follows that
\be\label{6.12}
\lim_{t\rightarrow 0}t^{\frac1q}\|u(t)\|_{L^p}=0.
\ee
Thus we arrive at

\begin{thm}\label{thm6.1}
Let $(q,p,r_0)$ be a generalized admissible triplet and let $\sigma={2\al}/b-n/p$.
Assume that $\varphi\in L^{r_0}(\mr^n).$
If $\ds\|\varphi\|_{\dot{B}^{-\sigma}_{p,\infty}}$ is small enough then the problem (\ref{3.2})
has a unique mild solution $u$ satisfying that
\ben
u(t,x)\in C\Big([0,\infty);L^{r_0}(\mr^n)\Big)\cap\mc{C}_q\Big([0,\infty);L^p(\mr^n)\Big).
\enn
Moreover, the solution $u$ satisfies (\ref{6.12}), that is,
\ben
u(t,x)\in C\Big([0,\infty);L^{r_0}(\mr^n)\Big)\cap\dot{\mc{C}}_q\Big([0,\infty);L^p(\mr^n)\Big).
\enn
\end{thm}

\begin{rem}\label{r6.1} {\rm
 Theorem \ref{thm6.1} indicates that the global solution
$u\in C\big([0,\infty);L^{r_0}(\mr^n)\big)$ exists provided
that the initial data $\varphi\in L^{r_0}(\mr^n)$ and its norm in the Besov space
$\ds\dot{B}^{-\sigma}_{p,\infty}(\mr^n)$ is small enough.
Note that the norm in $L^{r_0}(\mr^n)$ of $\varphi$ may be arbitrarily large.
For more details see the example in  \cite{Cannone}.}
\end{rem}

We may also consider the well-posedness in the Besov space $\ds\dot{B}^{-\sigma}_{p,\infty}(\mr^n)$
of the problem (\ref{3.2}). In doing so, we only need to use the solution space
\ben
X(I)=C_*([0,\infty);\dot{B}^{-\sigma} _{p,\infty}(\mr^n))\cap\mc{C}_q([0,\infty);L^p(\mr^n)
\enn
with the norm
\ben
\|u\|_{X(I)}=\sup_{t>0}\|u(t)\|_{\dot{B}^{-\sigma}_{p,\infty}}+\sup_{t>0}t^{\frac1q}\|u\|_p,
\enn
where $I=[0,\infty)$, $r_0={nb}/({2\al})\le p<r_0(b+1)$, $p>b+1$, $\sigma={2\al}/b-n/p\ge 0$
and $\ds\frac1q=\frac n{2\alpha}(\frac1{r_0}-\frac1p)$.
Since the Besov space $\ds\dot{B}^{-\sigma}_{p,\infty}(\mr^n)$ contains the self-similar initial
data $\varphi$, that is, $\ds\la^{\frac{2\alpha}b}\varphi(\lambda x)=\varphi(x)$ for any $\la>0$,
we also obtain self-similar solutions to the problem (\ref{3.2}) in this case.
By a similar argument as in the proof of Theorem \ref{thm6.1} we have the following
well-posedness in the homogeneous Besov space $\ds\dot{B}^{-\sigma}_{p,\infty}(\mr^n)$
of the problem (\ref{3.2}).

\begin{thm}\label{thm6.2}
Let $(q,p,r_0)$ be a generalized admissible triplet and let $p\ge r_0,$
$b+1<p<r_0(b+1)$. Assume that $\ds\varphi\in\dot{B}^{-\sigma}_{p,\infty}(\mr^n).$
Then, if $\ds\|\varphi\|_{\dot{B}^{-\sigma}_{p,\infty}}$ is sufficiently small then
the problem (\ref{3.2}) has a unique mild solution $u$ satisfying that
\ben
u\in C_*([0,\infty);\dot{B}^{-\sigma}_{p,\infty}(\mr^n))\cap\mc{C}_q([0,\infty);L^p(\mr^n).
\enn
Moreover, the solution $u$ satisfies (\ref{6.12}), that is,
\ben
u\in C_*([0,\infty);\dot{B}^{-\sigma}_{p,\infty}(\mr^n))\cap\dot{\mc{C}}_q([0,\infty);L^p(\mr^n).
\enn
\end{thm}

\begin{rem}\label{r6.2} {\rm
For the problem (\ref{4.6}) similar results to Theorems \ref{thm6.1} and \ref{thm6.2} hold
if we take $r_0={nb}/({2\alpha-1})$ and $\sigma=({2\alpha-1})/b-n/p$.
}
\end{rem}

\begin{rem}\label{r6.3} {\rm
(i) Consider the problem (\ref{1.1}) with $F(u)=-|u|^bu$, that is, the problem
\begin{eqnarray}\label{6.23}
\begin{cases}
u_t+(-\triangle)^\alpha u= -|u|^bu,\ (t,x)\in \mr^+\times\mr^n,\\
u(0,x)=\varphi(x),\ x\in \mr^n.
\end{cases}
\end{eqnarray}
It can be shown that the global solution to (\ref{6.23}) exists
under the condition $b<{4\al}/n$. In fact, multiplying both sides of the first equation of
(\ref{6.23}) by $u$ and integrating the equation thus obtained over $\mr^n$ gives
\ben
\frac12\frac d{\md t}\|u(t)\|^2_2+\|(-\tr)^{\frac{\alpha}2}u(t)\|^2_2+\|u(t)\|^{b+2}_{b+2}=0.
\enn
Thus we have
\ben
\|u\|_2\le C\|\varphi\|_2.
\enn
Therefore, let $(q,p,2)$ be an admissible triplet and define the solution space
$$
X(I)=C(I;L^2(\mr^n))\cap L^q(I;L^p(\mr^n))
$$
(or let $(q,p,2)$ be a generalized admissible triplet and define the solution space
$Y(I)=C(I;L^2(\mr^n))\cap\mc{C}_q(I;L^p(\mr^n))$), where $I=[0,T)$ for $T>0.$
Then by Theorem \ref{thm4.1} or \ref{thm4.2} the problem (\ref{6.23}) has a unique local solution
$u\in X(I)$ (or $u\in Y(I)$), where $T=T(\|\varphi\|_2)$ and $b<{4\al}/n$. Picard's method implies that
the local solution can be extended to be a global one.

(ii) If $\alpha=1$, then the restriction $b<{4\al}/n$ on the nonlinear growth can be removed.
}
\end{rem}

\section*{Acknowledgements} The research of Changxing
Miao was partly supported by the National Natural Science
Foundation (NNSF) of China. The research of Baoquan Yuan was
partly supported by the NNSF of China and the Science Foundation
for the Excellent Young Teachers of Henan Province. The research
of Bo Zhang was supported by the Chinese Academy of Sciences
through the Hundred Talents Program.

\end{document}